\documentclass[12pt]{amsart}
\usepackage{amsmath,amssymb,amscd,latexsym}

\def\ZZ{\mathbb Z}
\def\CC{\mathbb C}
\def\cO{\mathcal O}
\def\cS{\mathcal S}

\def\ra{\rightarrow}
\def\lra{\longrightarrow}
\def\PP{\mathbb P}
\def\QQ{\mathbb Q}
\def\op{\oplus}
\def\we{\wedge}
\def\ot{\otimes}
\def\a{\alpha}
\def\b{\beta}
\def\d{\delta}

\newtheorem{theorem}{Theorem}[section]
\newtheorem{proposition}[theorem]{Proposition}

\newtheorem{lemma}[theorem]{Lemma}
\newtheorem{corollary}[theorem]{Corollary}

\newtheorem{${}$}[theorem]{${}$}

\begin{document}

\title{Pfaffian lines and vector bundles \linebreak
       on Fano threefolds of genus 8}

\author{A. Iliev \ \and \ L. Manivel}
\date{}

\begin{abstract} Let $X$ be a general complex Fano threefold of genus 8. 
We prove that the moduli space of rank two semistable sheaves on 
$X$ with Chern numbers $c_1=1$, $c_2=6$ and $c_3=0$ is isomorphic to 
the Fano surface $F(X)$ of conics on $X$. This surface is smooth and 
isomorphic to the Fano surface of lines in the orthogonal to $X$ 
cubic threefold. Inside $F(X)$, the non-locally free sheaves are 
parameterized by a smooth curve of genus 26 isomorphic to the base 
of the family of lines on X. 

\end{abstract}

\maketitle

\footnotetext[0]{\ The 1-st author is partially supported by grant MI1503/2005
                  of the Bulgarian Foundation for Scientific Research}

\section{Introduction}\label{introduction}

A vector bundle $E$ on a smooth complex projective $n$-fold $X$ 
is {\it without intermediate cohomology} if $h^i(X,E(k)) = 0$ 
for any $i\ne 0,n$ and  $k\in\ZZ$. By a well known criterion 
of Horrocks, such a bundle must split if $X=\PP^n$. If $X=\QQ^n$
is a smooth quadric and $E$ is indecomposable, it must be a 
twist of a spinor bundle, see \cite{Kn,Ot}. 
In addition, as seen by Buchweitz, Greuel and Schreyer,
these are the only smooth $n$-folds with a finite number, 
up to twist, of indecomposable vector bundles without intermediate
cohomology. 

\smallskip
Until now more or less complete descriptions of these bundles 
have been obtained for only some restricted classes of varieties, 
see \cite{IM4} for a more comprehensive account of the known 
results. In particular, arithmetically Cohen-Macaulay (aCM) vector bundles
(that is, indecomposable rank two bundles with no intermediate
cohomology) on prime Fano threefolds have attracted considerable
attention; 
we refer to \cite{IP} for general facts about Fano varieties. 

\smallskip
Prime Fano threefolds $Y_d$ of index two and degree $d$ exist for 
$1 \le d \le 5$. The classification of aCM bundles on 
$Y_d$ is known for $d\ge 3$ \cite{AC}. Although the classification 
of aCM bundles on $Y_1$ (the double Veronese cone) and $Y_2$ (the 
quartic double solid) can easily be derived by similar methods,
it still remains unwritten.
The degree $d = -K_X^3 = 2g-2$ of a prime Fano threefold 
$X = X_d$ of index one is always even -- the integer 
$g$ is called the genus of $X$. Such Fano's exist 
for $2 \le g \le 12, g \not=11$, see e.g. \cite{Mu4}. 
The classification of all possible aCM bundles 
on these threefolds  (but not the existence of all of them) is given 
by C. Madonna in \cite{Ma} (see also \cite{AC}). There are finitely
many possible Chern numbers for a normalized aCM bundle on 
a prime Fano threefold $X$, and then a finite number of families 
for each possible choice, whose general
member is a stable vector bundle 
obtained by  Serre's construction from a subcanonical 
curve on $X$. Madonna deduces a list of 
$91$ possible pairs $(c_1,c_2)$, with $-1 \le c_1 \le 3$, corresponding to 
lines, conics and certain elliptic, canonical or half-canonical 
curves (provided such curves on $X_d$ exist).   

\smallskip
For a pair $(c_1,c_2)$ from the lists of Arrondo-Costa and Madonna,
denote by $M_X(2;c_1,c_2)$ the Maruyama moduli space of semistable 
rank two coherent sheaves on the Fano threefold $X$ with these 
Chern classes, and $c_3=0$. In some cases this moduli space happens to 
consist of a single point, corresponding to a unique bundle $E$. 
Such is the moduli space $M_X(2;1,5)$ for the Fano threefold $X=X_{14}$,
and the bundle $E$ defines a unique embedding of $X_{14}$ 
in the Grassmannian $G(2,6)$ parameterizing planes in a six dimensional
complex vector space. The study of rigid bundles on the Fano 
threefolds, K3 surfaces and canonical curves of genus $6,7,8,9,10$ is 
the base of Mukai's classifications by embeddings into homogeneous varieties.  

\smallskip
The next step is the study of moduli spaces $M_X(2;c_1,c_2)$ 
of non-rigid aCM bundles on prime Fano threefolds. 
The first attempts in this direction have been made in \cite{MT1} 
and \cite{IM1}, with the study of the moduli space $M_Y(2;0,2)$ 
on a cubic threefold $Y = Y_3$ based on a parallel study 
of the family of subcanonical curves -- elliptic quintics -- 
that appear as zero-loci of sections of the general $E \in M_Y(2;0,2)$. 
The Abel-Jacobi map sends the family of elliptic quintics on $Y$ onto the 
5-dimensional intermediate Jacobian $J(Y)$, and one can deduce 
that $M_Y(2;0,2)$ is birational to $J(Y)$. This statement was made
more precise by S. Druel \cite{D}, who proved that $M_Y(2;0,2)$ 
is isomorphic to the blow-up of $J(Y)$ along a copy of the Fano surface 
of lines in $Y$; see \cite{B1} for a joint presentation on both approaches. 

From some other point of view the same moduli space $M_Y(2;0,2)$ 
compactifies the set of Pfaffian representations of the 3-fold 
$Y$, an idea due originally to Adler, see \cite{AR}. 
As shown in \cite{IM2}, a similar question can be stated and answered 
similarly also for the quartic threefold $X = X_4 \subset\PP^4$, 
which by itself is a prime Fano 3-fold of index $1$: the set of 
Pfaffian representations of the 3-fold quartic $X$ is compactified by the 
7-dimensional moduli space of aCM bundles $M_X(2;3,14)$ on $X$.  
We refer to the paper \cite{B2} of Beauville for a modern and more 
general view on determinantal and Pfaffian representations of 
homogeneous forms. 

The ideas from \cite{MT1, IM2} have been used in 
\cite{IM3} and \cite{IM4} in the study of the moduli spaces $M_X(2;1,5)$ 
and $M_X(2;1,6)$ coming correspondingly from elliptic quintics and elliptic 
sextics on the prime Fano 3-fold $X = X_{12}$ of index $1$ and genus seven. 
To this circle of works one can include the new parameterization given by
Tikhomirov (see \cite{Tih}) of the theta divisor $\Theta$ for the quartic double 
solid $Y = Y_2$ by elliptic quintics, that in particular yields a birationality 
between $\Theta$ and a component of the moduli space $M_{Y}(2;0,3)$, 
see also \cite{MT2}. 
More recent, and from a different point of view, is the 
study by Arrondo and Faenzi of the aCM bundles on the prime Fano 3-folds $Y_5$ 
and $X_{12}$ in terms of monads, see \cite{AF} and \cite{F}. 

Note the common weak point in all these descriptions: they only 
consider the open subset of stable vector bundles in the moduli
space.

\smallskip
In this paper we give a full description of the moduli space 
$M_X(2;1,6)$ on the {\it general} Fano threefold $X = X_{14}$ 
of index one and genus eight. 
We combine the geometric approach used in 
\cite{IM1}-\cite{IM4} for other Fano threefolds with the ideas
of S. Druel \cite{D} to get not only all the 
vector bundles but also all the non locally free sheaves
parameterized by this moduli space. Our main result 
(Theorem \ref{main}) is that $M_X(2;1,6)$ is isomorphic to the
smooth surface $F(X)$ parameterizing the conics contained in $X$. 
Inside $F(X)$, the non locally free sheaves are parameterized 
by a smooth curve of genus 26 isomorphic to the family 
$\Gamma(X)$ of lines in $X$.  

Our proof is rather indirect, and we use the orthogonal cubic
threefold $Y$ of $X$ as an essential tool. Remember that $X$ is a 
generic linear section of the Grassmannian $G(2,6)\subset\PP^{14}$. 
The orthogonal cubic threefold $Y$ is then obtained as the 
section of the cubic Pfaffian hypersurface by the orthogonal 
linear subspace of the dual projective space. We prove that 
the surface $F(X)$ of conics in $X$ is isomorphic to the
Fano surface $F(Y)$ of lines in $Y$. Given a vector bundle in 
$M_X(2;1,6)$, we prove that it is generated by global sections,
and next we show how to construct a line $\ell\subset Y$ from an 
elliptic sextic obtained as the zero locus of a general section. 
Conversely, a given line $\ell \subset Y$ defines uniquely 
a codimension two singular linear section of $G(2,6)$, a special 
rational projection of which turns out to be isomorphic, 
somewhat unexpectedly, to the Grassmannian $G(2,5)$. 
Pulling-back the tautological rank two bundle on $G(2,5)$, we get a 
vector bundle $E_{\ell}\in M_X(2;1,6)$. Moreover, these two 
processes are inverse to each other. 

The Abel-Jacobi map allows to conclude the proof: $M_X(2;1,6)$ is 
mapped bijectively onto $F(X)$, which by itself is embedded in
the intermediate Jacobian $J(X)$ of $X$. Since both $M_X(2;1,6)$
and $F(X)$ are smooth surfaces, they are isomorphic. 

\medskip\noindent
{\it Acknowledgements}: We thank Claire Voisin and St\'ephane 
Druel for their help. 

\section{Preliminaries I : Some Grassmannian geometry}

\subsection{The Grassmannian $G(2,6)$} 
                           
Lots of important geometric properties of the Fano threefold $X_{14}$
come from the geometry of the Grassmannian $G(2,6)$. Here we state 
some of them that will be used later. 

We denote by $V$ a six-dimensional complex vector space, and by
$G(2,6)=G(2,V) \hookrightarrow \PP(\wedge^2V)=\PP^{14}$ 
the Grassmannian of planes in $V$, in its Pl\"ucker embedding. 
This is a smooth Fano manifold of dimension 8, degree $14$ and  
index 6. 

The secant variety of  $G(2,6)$ in $\PP^{14}$ is the cubic hypersurface 
of skew-symmetric tensors of rank at most four. Its equation is given 
by the Pfaffian. The action of $PGL_6$ in $\PP^{14}$ has only three
orbits, defined by the rank: $G(2,6)$, its complement in the Pfaffian 
hypersurface, and the complement of this hypersurface. 

\subsection{Lines and conics in $G(2,6)$}
A line in $G(2,6)$ is of the form $\PP(\ell\wedge L)$ for some
line $\ell$ in $V$ and some three dimensional subspace $L$ of 
$V$ containing $\ell$. The Fano variety of lines in $G(2,6)$
can thus be identified with the flag variety $F(1,3;6)$, with 
obvious notations. 

Similarly, the family of planes in $G(2,6)$ has two connected 
components, both homogeneous under the action of $PGL_6$. Of course
these planes contain lots of conics of $G(2,6)$. Apart from 
these, smooth conics are defined by two points of $G(2,6)$ plus two
concurring tangents. We deduce that we can find independent 
vectors $e_0,e_1,e_2,e_3\in V$ such that the span of our conic
is generated by $e_0\wedge e_1, e_2\wedge e_3, e_0\wedge e_2+
e_1\wedge e_3$. Similarly, a singular but reduced conic whose
linear span is not contained in the Grassmannian will be 
a union of two lines $\langle e_0\wedge e_1,e_0\wedge e_2\rangle$
and $\langle e_0\wedge e_1,e_1\wedge e_3\rangle$. The case of 
double lines is settled by the following lemma:

\begin{lemma}\label{double}
Let $L$ be a projective plane, meeting $G(2,6)$ along a 
double line. Then we can find four independent vectors 
$e_0,e_1,e_2,e_3$ in $V$ such that $L$ is generated
by $e_0\we e_1, e_0\we e_2, e_0\we e_3+e_1\we e_2$.
\end{lemma}

\proof Our hypothesis is that we have a line $\ell
=\langle e_0\wedge e_1,e_0\wedge e_2\rangle$, contained
in the plane $L\subset\PP^{14}$, such that the scheme 
intersection of $L$ with $G(2,6)$ is a double structure 
on $\ell$. 
Complete $e_0,e_1,e_2$ with three other vectors $e_3,e_4,e_5$
to get a basis of $V$. Then $L$ contains a unique tensor 
of the form $\omega =e_0\we u+\phi$ where $u\in\langle
e_3,e_4,e_5\rangle$ and $\phi$ does not involve $e_0$. 
For $\Omega=x_1e_0\wedge e_1+x_2e_0\wedge e_2
+y\omega$, the vanishing of 
$$\Omega\wedge\Omega=y(e_0\wedge (x_1e_1+x_2e_2+yu)\wedge\phi)
+y^2\phi\we\phi$$
defines our scheme-theoretic intersection. In order to get a 
double line, this equation must reduce to $y^2=0$, and 
we thus need that $e_1\wedge\phi=e_2\wedge\phi=0$. Hence 
$\phi$ must be a multiple of $e_1\we e_2$. Finally, for $L$ 
not to be contained in $G(2,6)$ we need $u$ not to be contained
in $\langle e_0,e_1,e_2\rangle$.\qed

\subsection{Lines in the Pfaffian hypersurface}
We will be interested in pencils of skew-symmetric tensors of 
constant rank four, i.e. lines in the Pfaffian hypersurface that 
do not meet the Grassmannian $G(2,6)$. 

We define an {\it A-line} to be a line generated by 
tensors of the form
\begin{eqnarray}
\nonumber & e_0\we e_2+e_1\we e_3, \\
\nonumber & e_0\we e_4+e_1\we e_5, 
\end{eqnarray}
for some basis $e_0,\ldots ,e_5$ of $V$.  
Similarly, a {\it B-line} will be generated by 
tensors of the form
\begin{eqnarray}
\nonumber & e_0\we e_2+e_1\we e_3, \\
\nonumber & e_0\we e_3+e_1\we e_4, 
\end{eqnarray}
for some independent vectors $e_0,\ldots ,e_4$ of $V$.

Both types of lines are contained in the tangent space to the
Grassmannian $G(2,6)$ at the unique plane $e_0\we e_1$. The 
main difference between the two types is that 
a B-line is contained, contrary to an A-line, 
in the projective span of a copy of $G(2,5)$ inside $G(2,6)$.  
As proved in \cite{manmez}, every line of
skew-symmetric tensors of constant rank four is an A-line 
or a B-line; moreover, 
the B-lines describe a hypersurface 
in the closure of the 22-dimensional family of A-lines.

For future use we note the following easy result, which can be 
obtained by an explicit computation. If $\ell$ is an A-line or 
a B-line, a point on $\ell$ represents a skew-symmetric form 
of rank four on $V^*$, whose kernel defines a projective line. 
Denote by $Q^{\ell}\subset\PP V^*$ the union of these lines. 

\begin{lemma}\label{notequal}
If $\ell$ is an A-line, then $Q^{\ell}$ is a smooth quadric surface 
in $\PP V^*$. 
If $\ell$ is a B-line, then $Q^{\ell}$ is a quadratic cone.
\end{lemma}

\section{Preliminaries II : Prime Fano threefolds of genus 8}

\subsection{Prime Fano threefolds from $G(2,6)$}
Since the Grassmannian $G(2,6)$ has dimension 8 and index 6, 
any smooth transverse linear section 
$$
X = X_{14} = G(2,6)\cap \PP^9_X \subset\PP^{14}
$$ 
is a Fano threefold of index one and degree $d = 14$, hence 
of genus $g = d/2 + 1 = 8$, i.e. the
smooth codimension $2$ linear sections of $X$  
are canonical curves of genus 8.
As shown independently by Gushel' and Mukai, 
see \cite{G}, \cite{Mu4}:  

\medskip

{\it Any smooth prime Fano threefold $X = X_{14}$ of index one 
and genus $8$ is a transverse section of 
$G(2,6) \subset \PP^{14}$ 
by a linear space $\PP^9_X \subset \PP^{14}$.}

\medskip
In Mukai's notation, the embedding $X \hookrightarrow G(2,6)$ is 
given by the unique rank 2 stable vector bundle $E_o\in M_X(2;1,5)$.
This bundle $E_o$ is defined by the Serre construction from any elliptic 
quintic on $X$, see \cite{Mu4}.
 
\medskip
The even Betti numbers of $X$ are $b_2=b_4=1$, canonical generators 
being given by the class of a hyperplane section and the class of a
line in $X$, respectively. The odd Betti numbers are $b_1=0$ and
$b_3=10$. In particular the intermediate Jacobian of $X$ is a five 
dimensional abelian variety $J(X)$. 

This can be seen as follows. Consider the 
normal exact sequence 
$$0\ra TX\ra TG_{X}\ra \cO_X(1)^{\oplus 5}\ra 0, $$
where $TG$ is the tangent bundle of the Grassmannian 
and $TG_{X}$ is its restriction to $X$. Note that $TX(-1)=\Omega_X^2$
since $K_X=\cO_X(-1)$. Twisting the previous sequence by 
$\cO_X(-1)$ and taking cohomology, we deduce an exact sequence 
$$H^0(TG_X(-1))\ra H^0(\cO_X^{\oplus 5})\ra H^1(\Omega_X^2)\ra
H^1(TG_X(-1)).$$
Since $X$ is a linear section of the Grassmannian there is an
associated Koszul complex that we can use to deduce that 
$H^i(TG_X(-1))=0$ whenever $H^{i+j}(TG(-j-1))=0$ for $0\le j\le 5$. 
For $i=0$ and $i=1$ this easily follows from Bott's theorem (see e.g. 
\cite{man}, Proposition 2). 
We conclude that $h^{1,2}(X)=5$, as claimed.

\subsection{The orthogonal cubic threefold}
For the Fano threefold $X=G(2,6) \cap \PP_X^9$,  
consider the orthogonal four dimensional space 
in the dual projective space $\PP(\wedge^2V^*)$. 
This four dimensional projective space meets the 
Pfaffian hypersurface in $\PP(\wedge^2V^*)$ along a 
cubic threefold $Y\subset\PP^4_Y$, the {\it orthogonal cubic threefold} 
of the Fano threefold $X$. 

For $X$ general, $\PP^4_Y$ is a general subspace in $\PP^{14}$,
hence does not meet the Grassmannian $G(2,V^*)$, whose codimension 
is 6. Therefore any point $y\in Y$ represents a skew-symmetric form 
of rank four on $V$, whose kernel is a well-defined projective
line $n_y\subset\PP^5$.  

Since the representation of $X$ as a linear section 
of $G(2,6)$ is unique up to the action of $PGL_6$ on 
the codimension five spaces in $\PP^{14}$, the dual 
cubic threefold $Y$ is uniquely defined up to the action 
of $PGL_6$. Moreover, $Y$ is smooth whenever $X$ is smooth, 
see \cite{Pu}.

%Equivalently, we get a rank two vector bundle $\Sigma_X$ on $Y$, 
%which belongs to the moduli space $M_Y(2;0,2)$ of semistable
%rank two sheaves on $Y$ with $c_1=0, c_2=2, c_3=0$.  
%As shown by Druel, this moduli space is isomorphic to 
%a blow-up of the intermediate Jacobian $J(Y)$ along  (minus)
%the translate 
%of the Abel-Jacobi image of the Fano surface $F(Y)$ \cite{D}.   
%The vector bundles in this moduli space  
%parameterize the Pfaffian representations of the
%cubic threefold $Y$, see Appendix V in \cite{AR} and 
%Prop.8.5 in \cite{B2}.

\subsection{The Palatini quartic}
For $x\in X$, denote by $\ell_x \subset\PP^5$
the corresponding line. Let  
\begin{eqnarray}
\nonumber 
 & W = &\cup \{ \ell_x : x \in X \} \subset \PP^5, \\ 
\nonumber 
 & V = &\cup \{ n_y :  y \in Y \} \subset \PP^5.
\end{eqnarray}
Both $W$ and $V$ are subvarieties of $\PP^5$
swept out by lines, and in fact they coincide. 
More precisely, the following takes place
(see \S 50 in \cite{AR} or \cite{Pu}):   
%P. J. Puts: page 83

\medskip

%\begin{lemma}\label{gen-3} 
{\it 
{\bf (i)}
$W = V$ is an irreducible quartic hypersurface in $\PP^5$, 
whose singular locus is a curve $\Gamma(W)$ of degree $25$; 

{\bf (ii)}
through any point $p \in W - \Gamma(W)$  passes a unique
line $\ell_x, x \in X$ and a unique kernel line $n_y, y \in Y$;  
and through any point $v \in \Gamma(W)$ passes exactly a 
linear pencil of lines $\ell_x, x \in X$. 
}
%\end{lemma}

%\noindent
%{\bf Proof}: p.83 in \cite{Pu} and p.174 in \cite{AR}.  

\medskip
The quartic hypersurface $W \subset \PP^5$ is known as
the {\it Palatini quartic} of the prime Fano threefold $X$ of genus $8$. 
The assertion (ii) above implies that a general hyperplane section 
of $W$ is a singular quartic threefold with $25$ singular points  
which is birational to both $X$ and $Y$. 
In particular $X$ and $Y$ are birational (and both unirational but 
not rational). 

This has an interesting consequence, observed by Puts, 
that will be of crucial use later.

\begin{proposition}\label{isojac}
The intermediate Jacobians of $X$ and $Y$ are isomorphic. 
\end{proposition}

\proof We briefly recall the proof given in \cite{Pu}.
Let $H$ be a general hyperplane in $\PP^5$. Define the 
rational map $u_{\scriptscriptstyle X} : X\dasharrow W\cap H$ 
by mapping $x\in X$ to the intersection point of the line 
$\ell_x$ with $H$. Similarly 
$u_{\scriptscriptstyle Y} : Y\dasharrow W\cap H$
maps $y\in Y$ to the intersection point of the kernel line 
$n_y$ with $H$. Both maps are birational, as we just noticed. 

The map $u_{\scriptscriptstyle X}$ is well defined 
at $x$ if the line $\ell_x$ is not contained in $H$, 
that is,  outside the intersection $C=X\cap G(2,H)$. 
Since $H$ is general this is a smooth elliptic quintic. 
The map $u_{\scriptscriptstyle Y}$ is also defined 
outside a smooth elliptic quintic $D\subset Y$, 
and there exists a commutative diagram 
$$\begin{array}{ccccc}
X & \stackrel{{\scriptstyle u}_{\scriptscriptstyle X}}{\rightarrow} & W\cap H & 
\stackrel{{\scriptstyle u}_{\scriptscriptstyle Y}}{\leftarrow} & Y \\
\uparrow & & & & \uparrow \\
X' & \leftarrow & Z & \rightarrow & Y' 
\end{array}$$
where the vertical maps are the blow-up of the quintic elliptic
curves $C$ and $D$, and the arrows $Z\ra X'$ and $Z\ra Y'$ are the
blow-up of 25 lines. By \cite{CG}, one deduces that 
$$J(X)+J(C)\simeq J(X')\simeq J(Z)\simeq J(Y')\simeq J(Y)+J(D).$$
This implies that $J(C)\simeq J(D)$, hence $C$ and $D$ are isomorphic 
(see also Ch.III \S 1 in \cite{Is}), 
and also that our birationalities induce an isomorphism 
$J(X)\simeq J(Y)$. \qed

\section{Lines on $X_{14}$ and its orthogonal cubic}

\subsection{Lines on the orthogonal cubic threefold}\label{F(Y)} 
 
We recall some known facts about the family
of lines on the $3$-dimensional cubic hypersurface. 
The basic reference about cubic threefolds and intermediate 
Jacobians of threefolds is the paper 
\cite{CG} of Clemens and Griffiths.

\smallskip
{\it 
Let $Y \subset \PP^4$ be a general cubic $3$-fold. Then

\smallskip
{\bf (i)}  
The family $F(Y)$ of lines on $Y$
is a smooth irreducible surface of general type, 
of geometric genus $p_g = 10$, of irregularity 
$q = 5$, and such that $K^2 = 45$. The canonical system $|K|$ 
is very ample and defines the Pl\"ucker embedding 
$F(Y) \subset G(2,5) \rightarrow \PP^9$. 

\smallskip
{\bf (ii)}   
The intermediate Jacobian $J(Y) = H^{2,1}(Y)^*/H_3(Y,\ZZ)$ 
is an abelian variety of dimension $5$, and the Abel-Jacobi map 
$\Phi : F(Y) \rightarrow J(Y)$ is an embedding. In particular, 
$F(Y)$ contains no rational curve.}

\smallskip
When $Y$ is the orthogonal cubic threefold to a Fano threefold
$X$ of genus 8, a line in $Y$ is a pencil of skew-symmetric
forms on $V^*$, all of rank four. Since the B-lines  
form a codimension one family of lines of skew-symmetric
forms of constant rank four, we get:  

\smallskip
{\it {\bf (iii)} 
The family of B-lines in the orthogonal cubic threefold $Y$
to a general Fano threefold $X$ of genus 8 is a smooth curve
$\Gamma(Y)\subset F(Y)$.}  

\smallskip
%\noindent {\bf Remark}. 
For the smoothness of $\Gamma(Y)$, see \cite{AR}, Corollary (49.8)
and Lemmas (51.3) and (51.15). 
Notice that the curve $\Gamma(Y)$ is not invariantly defined 
by the cubic threefold $Y$; as a subset of $F(Y)$ 
it depends on the choice of $X$ in the 5-dimensional family 
of Fano threefolds whose orthogonal cubic is $Y$. 
%As shown in \cite{AR}, Appendix V, this family of 
%Fano threefolds  can be identified with  the set of 
%Pfaffian representations of the cubic $Y$.
%
Other characterizations of B-lines are given in 
\cite{AR}. They are the lines $\ell\subset Y$ such that 
the kernel lines  $n_y, y \in\ell$  have a common point
(see Definition (49.1) in \cite{AR}). Equivalently, they 
are the {\it jumping lines} of the restriction to $Y$ of 
the tautological rank two bundle (see Lemma (49.2) 
in \cite{AR}).

\subsection{Lines in the general Fano threefold $X_{14}$}\label{xlines}

It is well known that the scheme $\Gamma(X)$ 
of lines on the general Fano threefold $X$ of 
genus 8 is a smooth irreducible curve of genus 26, 
see e.g. Prop. 4.2.2 and Th. 4.2.7 in \cite{IP};    
in particular the normal bundle to any line $\ell\subset X$ 
is $N_{\ell/X}=\cO_{\ell}\op\cO_{\ell}(-1)$, ibid. 
This implies that $h^0(N_{\ell/X})=ext^1(\cO_{\ell},\cO_{\ell})=1$
and $h^1(N_{\ell/X})=ext^2(\cO_{\ell},\cO_{\ell})=0$.

\smallskip
For a line $\ell$ in $G(2,6)$, denote by $v(\ell)\in\PP^5$ the 
{\it vertex} of $\ell$, the common point of the lines in $\PP^5$
parameterized by $\ell$. 
Let $H(\ell)\subset V^*$ denote the hyperplane orthogonal
  to the vertex $v(\ell)$. 
Then $\PP(\wedge^2H(\ell))\simeq \PP^9$ has codimension three
in $\ell^{\perp}$, which also contains $\PP^4_Y$. Therefore these
two subspaces of $\ell^{\perp}$ have to meet
along some projective space of dimension at least one. Moreover
this space is contained in $Y$ since $H(\ell)$ has dimension five, 
and in dimension five a skew-symmetric form has rank at most four. 
But $Y$ does not contain any plane, so 
$\PP(\wedge^2H(\ell))\cap \PP^4_Y$ is a projective line $\ell '$ in $Y$, 
which is obviously a B-line. 

Conversely, let $\ell' \subset Y$ be a B-line, hence contained in the
span $\PP(\wedge^2H)$ of $G(2,H)$, a copy of $G(2,5)$, for some 
hyperplane $H$ of $V^*$. The orthogonal to $H$ is point $v\in\PP^5$, 
and the Schubert cycle of lines passing through $v$ is a $\PP^4_v$
contained in $G(2,6)$. Both $\PP^4_v$ and $\PP^9_X$ are contained in 
the 12-space $(\ell')^{\perp}$, so they must meet each other 
along a linear space $\ell$ of dimension at least one. Since $\ell$ 
is contained in $X = G(2,6) \cap \PP^9_X$, and since $X$ contains no planes,
then $\ell$ must be a line. 
Together with the preceding, this yields the following:

\begin{proposition}\label{xlines-blines}
The curve $\Gamma(X)$ of lines in $X$ is isomorphic to 
the curve $\Gamma(Y) \subset F(Y)$ of B-lines  
in the orthogonal cubic $Y$ of $X$.
\end{proposition}

Notice that if $\ell$ is a line in $X$, 
its vertex $v(\ell)$ is a singular point of the Palatini quartic $W$,
see Lemma (51.5) in App.V of \cite{AR}. 
Moreover the map $\Gamma(X)\ra\Gamma(W)$ sending a line
to its vertex is injective, since two lines with the same 
vertex would generate a plane contained in $X$. This implies:

\begin{proposition}\label{xlines-wsing}
The curve $\Gamma(X)$ of lines in $X$ is isomorphic to 
the singular locus $\Gamma(W)$ of the Palatini quartic of $X$. 
\end{proposition}

The fact that for $X$ general, $\Gamma(W)$ is a smooth curve
of genus 26 was already observed by Puts \cite{Pu}.

\subsection{Conics on $X_{14}$}\label{conics}

We denote by $F(X)$ the Hilbert scheme of conics on $X$, 
that is, the Hilbert scheme of closed 
subschemes of  $X$ with  Hilbert polynomial $P(n)=2n+1$.
It is known that $F(X)$ is reduced of pure dimension two:
it is called traditionally the {\it Fano surface} of $X$. 
Any general point of $F(X)$ parametrizes a smooth conic 
$q \subset X$ with normal bundle 
$N_{q|X} = {\mathcal O}_q \oplus {\mathcal O}_q$. 
Moreover the closed subset of singular conics 
(pairs of intersecting lines, or double lines) is of pure dimension $1$.
See \S 4.2 in \cite{IP}.  

Of course any conic on $X$ is a conic on $G(2,6)$, and 
the different types of conics in $G(2,6)$ can easily
be described. We can distinguish three types:
\begin{enumerate}
\item Conics whose linear span is a plane contained in $G(2,6)$;
\item Reduced plane sections of a copy of $G(2,4)$ in $G(2,6)$;
\item Double lines.
\end{enumerate}
Clearly conics of the first type cannot be contained in a general
Fano threefold $X$, which contains no planes. 
Double lines are handled as follows:

\begin{lemma}\label{nodouble}
The general Fano threefold $X_{14}$ contains no double lines.
\end{lemma}

\proof Double lines on $G(2,6)$ were described in Lemma \ref{double}. 
An immediate consequence is that there is a $\PP^3$ of projective 
planes whose scheme intersection with $G(2,6)$ is a double line 
supported on a given line in $G(2,6)$. Since the Fano variety of 
lines on the Grassmannian is 11-dimensional, we get a family
of planes of dimension 14, and the family of $\PP^9$'s containing
a plane in this family has dimension $14+7\times 5=49$. This is one
less than the dimension of the Grassmann variety of $\PP^9$'s in 
$\PP^{14}$, so our claim follows. \qed

\begin{proposition}\label{F(X)-smooth} 
The Fano surface $F(X)$ of conics on the general $X$ 
is smooth.
 \end{proposition}

\proof
Let $q\subset X$ be a reduced conic. To prove that 
the $F(X)$ is smooth and two-dimensional
at $q$, we will check that $h^1(N_{q/X})=0$ and $h^0(N_{q/X})=2$, where 
$N_{q/X}$ denotes the rank two normal bundle of $q$ in $X$
(note that $q$ can be singular but is always a locally complete
intersection). 

\medskip
{\it First case : $q$ is smooth}. 
We know that $q=G(q)\cap X$,
where $G(q)$, the Grassmannian of lines in $\PP^3_q$, 
is a four dimensional quadric
in $G=G(2,6)$. We have an exact sequence
$$0\ra N_{q/G(q)}\ra N_{q/G}\ra N_{G(q)/G|q}\ra 0.$$
Since $q$ is a transverse linear section of $G(q)$, we have
$N_{q/G(q)}=\cO(2)^3$. The restriction of the tautological 
bundle $T^*$ to $q$ is globally generated and has degree two, 
thus
it must decompose into $\cO\op\cO(2)$ or $\cO(1)\op\cO(1)$. 
In the first case, the lines parametrized by $q$ would have a common
point, and the linear span of $q$ would be contained in $G$, 
hence in $X$, which is impossible. So we must be in the second case, and
since $G(q)$ is the zero
locus of a section of $T^*\oplus T^*$ on $G$, we deduce that 
$N_{G(q)/G|q}=\cO(1)^4$.

\smallskip
Now, consider the commutative diagram:
$$\begin{array}{ccccccccc}
& & 0 & & 0 & & & & \\
& & \downarrow & & \downarrow & & & & \\
0 &\ra & Tq & \ra &TG(q)_{|q} &\ra  & N_{q/G(q)}& \ra & 0 \\
& & \downarrow & & \downarrow & &\downarrow & & \\
0 &\ra & TX_{|q} & \ra &TG_{|q} &\ra  & N_{X/G|q} & \ra & 0 \\
& & \downarrow & & \downarrow & & & & \\
& & N_{q/X} & \ra & N_{G(q)/G|q} & & & & \\
& & \downarrow & & \downarrow & & & & \\
& & 0 & & 0 & & & & 
\end{array}$$
Suppose that the image $\bar{z}\in N_{q/X}$ of some $z\in TX_{|q}$,
is mapped to zero in $N_{G(q)/G|q}$. From the diagram, we deduce
that $z\in TG(q)_{|q}$, hence $z\in Tq=TX_{|q}\cap TG(q)_{|q}$. But
this means that $\bar{z}=0$. We conclude that the map 
$N_{q/X}\ra N_{G(q)/G|q}=\cO(1)^4$ is injective. Since $N_{q/X}$
has degree zero, this implies that $N_{q/X}=\cO\op\cO$ or 
$N_{q/X}=\cO(1)\op\cO(-1)$, hence our claim. 

\medskip
{\it Second case: $q$ is singular}. 
So $q$ is the union of two coplanar
but distinct lines $\ell$ and $m$, meeting at a point $p$. 
We prove that $Ext^2(I_q,I_q)=0$ 
and $Ext^1(I_q,I_q)$ is two dimensional. 
We begin with the short exact sequence 
$0\ra I_q\ra I_{\ell}\ra\cO_m(-p)\ra 0$. Applying
$Hom(I_q,.)$, we get the long exact sequence
\begin{eqnarray}\nonumber
0\ra Hom(I_q,I_q)\ra Hom(I_q,I_{\ell})\ra
Hom(I_q,\cO_m(-p))\ra \\ 
\hspace{1cm}\ra Ext^1(I_q,I_q)
\ra Ext^1(I_q,I_{\ell})\ra Ext^1(I_q,\cO_m(-p))\ra  \\
\nonumber \hspace{2cm}\ra Ext^2(I_q,I_q)\ra
Ext^2(I_q,I_{\ell})\ra Ext^2(I_q,\cO_m(-p)).
\end{eqnarray}
Applying $Hom(.,\cO_m(-p))$, we obtain
\begin{eqnarray}\nonumber 
0\ra Hom(\cO_m,\cO_m)\ra 
Hom(I_{\ell},\cO_m(-p))\ra Hom(I_q,\cO_m(-p))\ra \\
 \ra Ext^1(\cO_m,\cO_m)\ra 
Ext^1(I_{\ell},\cO_m(-p))\ra Ext^1(I_q,\cO_m(-p))\ra \\
\nonumber
\ra Ext^2(\cO_m,\cO_m)\ra 
Ext^2(I_{\ell},\cO_m(-p))\ra Ext^2(I_q,\cO_m(-p)).
\end{eqnarray}
Note that $Hom(\cO_m,\cO_m)=Hom(I_{\ell},\cO_m(-p))
=\CC$, while $Hom(I_q,\cO_m(-p))=0$. 

To compute $Ext^k(\cO_m(-p),\cO_m(-p))=Ext^k(\cO_m,\cO_m)$, 
we use the short exact sequence $0\ra I_m\ra \cO_X\ra\cO_m\ra 0$. Since 
$Ext^k(\cO_X,I_m)=H^k(I_m)=0$ for $k\ge 0$ we deduce that 
$Ext^k(\cO_m,I_m)=Ext^{k-1}(I_m,I_m)$. Using 4.2,
we get that this group is equal
to $\CC$ for $k=1,2$ and to zero for $k=0,3$. 

Then we apply $Hom(\cO_m,.)$ to the previous short sequence:
\begin{eqnarray}\nonumber 
0\ra Hom(\cO_m,I_m)\ra 
Hom(\cO_m,\cO_X)\ra Hom(\cO_m,\cO_m)\ra \\
\ra Ext^1(\cO_m,I_m)\ra 
Ext^1(\cO_m,\cO_X)\ra Ext^1(\cO_m,\cO_m)\ra \\
\nonumber
\ra Ext^2(\cO_m,I_m)\ra 
Ext^2(\cO_m,\cO_X)\ra Ext^2(\cO_m,\cO_m).\;\;\;
\end{eqnarray}
We deduce that $Ext^k(\cO_m,\cO_m)=\CC$ for $k=0,1$, and 
zero for $k=2$. 

Now we apply the functor $Hom(.,\cO_m(-p))$ to the 
sequence $0\ra I_{\ell}\ra \cO_X\ra\cO_{\ell}\ra 0$.
Since $Ext^k(\cO_X,\cO_m(-p))=H^k(\PP^1,\cO(-1))=0$ for all $k$, 
we deduce that $Ext^k(I_{\ell},\cO_m(-p))=
Ext^{k+1}(\cO_{\ell},\cO_m(-p))$. To compute the latter, 
we use the existence of  a spectral sequence 
$$H^i(X,{\mathcal Ext}^j(\cO_{\ell},\cO_m(-p)))
\Rightarrow  Ext^{i+j}(\cO_{\ell},\cO_m(-p)).$$
Since the sheaf ${\mathcal Ext}^j(\cO_{\ell},\cO_m(-p))$ is 
supported at $p$, we deduce that 
$$Ext^{k+1}(\cO_{\ell},\cO_m(-p))=H^0(X,
{\mathcal Ext}^{k+1}(\cO_{\ell},\cO_m(-p))).$$
So we just need to compute the rank of the fiber of this local 
${\mathcal Ext}$-sheaf at the point $p=\ell\cap m$. Locally around
that point, $\ell$ and $m$ are given by some equations $x=y=0$
and $x=z=0$ respectively, and $\cO_{\ell}$ has a free resolution 
given by the Koszul complex defined by $x=z=0$. Applying 
${\mathcal Hom}(.,\cO_m(-p)))$ we get a complex whose 
cohomology is readily computed. We deduce that the rank of 
${\mathcal Ext}^{k+1}(\cO_{\ell},\cO_m(-p)))$ at $p$ is one 
for $k=1,2$, and zero for $k=3$. Thus $Ext^k(I_{\ell},\cO_m(-p))=\CC$
for $k=0,1$ and zero for $k=2$.

From the long exact sequence (2), we deduce that 
$Ext^k(I_q,\cO_m(-p))=0$ for $k=1,2$, and consequently from (1)
that  $Ext^k(I_q,I_q)=Ext^k(I_q,I_{\ell})$ again for $k=1,2$. 
To compute this group, we apply the 
functor $Hom(., I_{\ell})$ to the short sequence 
$0\ra I_q\ra I_{\ell}\ra\cO_m(-p)\ra 0$. This gives
\begin{eqnarray}\nonumber 
0\ra Hom(\cO_m(-p),I_{\ell})\ra 
Hom(I_{\ell},I_{\ell})\ra Hom(I_q,I_{\ell})\ra \\
\ra Ext^1(\cO_m(-p),I_{\ell})\ra 
Ext^1(I_{\ell},I_{\ell})\ra Ext^1(I_q,I_{\ell})\ra \\
\nonumber
\ra Ext^2(\cO_m(-p),I_{\ell})\ra 
Ext^2(I_{\ell},I_{\ell})\ra Ext^2(I_q,I_{\ell})\ra 0.
\end{eqnarray}
The last arrow is surjective because 
$Ext^3(\cO_m(-p),I_{\ell})$, since $m$ is a locally 
complete intersection of codimension two. Since ${\ell}$
is supposed to define a smooth point of the family of 
lines on $X$, we know that $Ext^2(I_{\ell},I_{\ell})=0$
and $Ext^1(I_{\ell},I_{\ell})=\CC$. Finally, we apply 
$Hom(\cO_m(-p),.)$ to $0\ra I_{\ell}\ra \cO_X\ra\cO_{\ell}\ra 0$
to obtain
\begin{eqnarray}\nonumber 
0\ra Hom(\cO_m(-p),I_{\ell})\ra 
Hom(\cO_m(-p),\cO_X)\ra Hom(\cO_m(-p),\cO_{\ell})\ra \\
\ra Ext^1(\cO_m(-p),I_{\ell})\ra 
Ext^1(\cO_m(-p),\cO_X)\ra Ext^1(\cO_m(-p),\cO_{\ell})\ra \\
\nonumber
\ra Ext^2(\cO_m(-p),I_{\ell})\ra 
Ext^2(\cO_m(-p),\cO_X)\ra Ext^2(\cO_m(-p),\cO_{\ell})\ra 0.
\end{eqnarray}
By Serre duality, $Ext^k(\cO_m(-p),\cO_X)=Ext^k(\cO_m(-2p),\omega_X)$ is
dual to the cohomology group 
$H^{3-k}(X,\cO_m(-2p))=H^{3-k}(\PP^1,\cO(-2))$, hence equal to
$\CC$ for $k=2$ 
and zero otherwise. We deduce that $Ext^1(\cO_m(-p),I_{\ell})=0$ and 
$Ext^2(\cO_m(-p),I_{\ell})=\CC$. Using (4) we get that 
$Ext^2(I_q,I_{\ell})=Ext^2(I_q,I_q)=0$, while 
$Ext^1(I_q,I_{\ell})=Ext^1(I_q,I_q)$ is two-dimensional. 
This concludes the proof that $F(X)$ is smooth at $q$. \qed

\medskip
 Let $q\in F(X)$ be a conic in $X$, possibly singular. There
 is a unique codimension two subspace $L$ of $V$ such that
 $q$ is a linear section of $G(2,L)$, and the linear space 
 $\PP(\wedge^2L)\simeq \PP^5$ meets $\PP_X^9$ along the
 span $\PP^2_q$ of the conic. Dually, this implies that the 
 orthogonal to $\PP(\wedge^2L)$, which is the linear space 
 $\PP(L^{\perp}\wedge V^*)\simeq \PP^8$, meets $\PP^4_Y$ along 
 a line $\ell$. But a skew-symmetric form in $L^{\perp}\wedge V^*$
 has rank at most four, so the line $\ell$ is in fact contained 
 in $Y$. 

 Conversely, a line $\ell$ in $Y$, of type A or B, is contained 
 in the tangent space $\PP(M\wedge V^*)$ of the Grassmannian 
 at a unique line $m=\PP M$. Dually, the orthogonal $\PP(\wedge^2
 M^{\perp})$ to this tangent space is contained in $\ell^{\perp}$. 
 The linear span $\PP^9_X$ of $X$, which has codimension three
 in $\ell^{\perp}$, has to cut the four dimensional quadric 
 $G(2,M^{\perp})\subset\PP(\wedge^2 M^{\perp})$ at least along a 
 conic. But since $X$ contains no surfaces of degree two, the 
 intersection $G(2,M^{\perp})\cap\PP^9_X=q$ is nothing 
 more than a conic. We have proved:

 \begin{proposition}\label{conics-5}
 The Fano surface $F(X)$
 of conics in $X$ is isomorphic with the 
 Fano surface $F(Y)$ of lines in the orthogonal cubic $Y$.
 \end{proposition}

Recall that the Albanese variety of $F(Y)$ is isomorphic
to the intermediate Jacobian $J(Y)$ of $Y$ \cite{CG}. 
Moreover, by Proposition \ref{isojac}, $J(Y)$ is isomorphic to the 
intermediate Jacobian $J(X)$ of $X$. Since $F(X)$ is smooth, 
the Abel-Jacobi map $F(X)\ra J(X)$ is algebraic. 

\smallskip
How is the Abel-Jacobi image of $F(X)$ related with $F(Y)$
in $J(X)\simeq J(Y)$?

Let $\ell$ be a line on $Y$. We have seen in Lemma \ref{notequal} that 
the union of the kernels $n_y$ of the points $y\in\ell$ is a 
two-dimensional quadric $Q^{\ell}$ of rank at least three, contained in 
the Palatini quartic $W$. The linear span of $Q^{\ell}$ is a 
codimension two linear space $\PP^3_{\ell}\subset\PP^5$.

On the side of $X$, we have a conic $q_{\ell}$
parametrizing lines in $\PP^5$ whose union is again a  
quadric $Q_{\ell}$, also contained in $W\cap\PP^3_{\ell}$. Note that 
when the conic $q_{\ell}$ is a union of two lines, $Q_{\ell}$ is a union 
of two planes. In particular $Q_{\ell}\neq Q^{\ell}$, and this must 
remain true for any general $\ell$. 

Note moreover that $W$ cannot contain $\PP^3_{\ell}$, since otherwise 
$Sing(W)\cap\PP^3_{\ell}$ would contain a surface or a curve of 
degree nine, in contradiction with the fact that  $\Gamma(W)=Sing(W)$
is a smooth curve of degree $25$. 

Since $W$ has degree four, we deduce that 
$$W\cap\PP^3_{\ell}=Q_{\ell}+Q^{\ell}$$
for all $\ell\in F(Y)$. 
Let $H\subset\PP^5$ be the general hyperplane that we used 
in Proposition \ref{isojac} to define
a birationality between $X$ and $Y$, and then an isomorphism between
$J(X)$ and $J(Y)$. The image of $\ell$ by $u_Y$ is the plane conic 
$Q^{\ell}\cap H$. The image of $q_{\ell}$ by $u_X$ is the plane 
conic $Q_{\ell}\cap H$. Since the sum of these two cycles is a 
linear section of $W\cap H$, the Abel-Jacobi image of $\ell$ in $J(Y)$
is equal to minus the Abel-Jacobi image of $q_{\ell}$ in $J(X)$, 
once we have identified 
the two intermediate Jacobians via Proposition \ref{isojac}.

Since the Abel-Jacobi map embeds $F(Y)$ in $J(Y)$, and since $F(X)$
and $F(Y)$ are isomorphic by Proposition \ref{conics-5},  we deduce:

\begin{proposition}\label{-1}
The Abel-Jacobi map from $F(X)$ to $J(X)$ is an embedding. 
Moreover,
the images of $F(X)$ in $J(X)$ and $F(Y)$ in $J(Y)$ are the 
same up to multiplication by $-1$. 
\end{proposition}

 \section{The moduli space $M_X(2;1,6)$} 
 %             and elliptic sextics on $X_{14}$}\label{mod} 

In this section $X = X_{14} = G(2,6) \cap \PP^9_X$ 
is again a general smooth prime Fano threefold of genus 8,
and we study the moduli space $M_X(2;1,6)$ of semistable rank two 
torsion free sheaves $E$ on $X$ with Chern numbers $c_1=1$, 
$c_2=6$ and $c_3=0$. In particular $c_1(E) = h$ is the hyperplane 
class of $X$, and since $h$ is not divisible in $Pic(X) = \ZZ h$,   
any such sheaf $E$ is in fact stable, and the stability 
of $E$ is equivalent to its slope stability.

 \subsection{Vector bundles in $M_X(2;1,6)$}

 First we study the locally free sheaves $E$ from
 the moduli space $M_X(2;1,6)$. 
 Note that by the stability assumption the bundle 
 $E^* = E(-1) = Hom(\cO_X(1),E)$ has no sections. 
 Moreover, the Riemann-Roch formula for vector bundles on threefolds 
 gives $\chi(E) = 5$,  
 see Ch.15.2: Example 15.2.5 in \cite{Fu}. 

 \begin{lemma}\label{mod-5}
 $h^2(E)=0$ and $h^0(E)>0$. 
 \end{lemma}

 \proof 
 Since $E(-2)=E^*(-1)$ has no sections, $h^3(E)=0$ by Serre duality. 
 Suppose that $h^2(E)=h^1(E^*(-1))\neq 0$. Then we would have a non
 trivial extension
 $$
 0 \rightarrow {\mathcal O}_X(-1) \rightarrow F \rightarrow E \rightarrow 0,
 $$
 where the rank three vector bundle $F$ 
 has $c_1(F) = 0,\ c_2(F) = -8$ and $c_3(F) = 6$. 
 Thus $\Delta(F) = deg(2c_1(F)^2 - 6c_2(F)) = 48 > 0$, 
 so by Bogomolov's inequality \cite{Bo}, $F$ cannot be semistable.
 So either $\cO_X(1)$ maps non trivially to $F$, hence to $E$ -- a
 contradiction!, or $F$ maps non trivially to $\cO_X(-1)$ and the
 extension splits, again a contradiction. 
 Thus $h^2(E)=0 \mbox{ and } h^0(E)=5+h^1(E)>0$. \qed

 \medskip
 Now we choose a general hyperplane section $S$ of $X$, and denote by 
 $E_S$ the restriction of $E$ to $S$. The Picard group of this K3 surface
 is generated by $\cO_S(1)$, and by Maruyama's restriction theorem $E_S$ 
 remains semistable, hence stable, with Chern classes 
 $c_1(E_S) = h_S$, the hyperplane class of $S$, and $c_2(E_S) = 6$, 
see \cite{Mar}. 

 \begin{proposition}\label{mod-6}
 \hspace*{2cm}
 \begin{enumerate}
 \item $h^0(E_S) = 5$, and $h^i(E_S) = 0$ for $i > 0$.
 %{\bf (ii)} 
 \item $h^0(E)=5$ and $h^1(E)=h^1(E(-1))=0$.
 \end{enumerate}
 \end{proposition}

 \proof (1)
 By Riemann-Roch for vector bundles on surfaces 
 $\chi(E_S)=5$ (see Example 15.2.2 in \cite{Fu}); 
 and since $H^2(E_S)$ and $H^0(E_S(-1))$ vanish by the 
 stability of $E_S$ then $h^0(E_S) = 5 + h^1(E_S) \ge 5.$ 

 In particular $E_S$ has non-trivial sections; 
 let $s$ be one of them, 
 and let $Z = Z(s)$ be the zero-scheme of $s$. The
 stability of $E_S$ implies that 
 $Z$ is either of codimension $2$ in $S$ or empty.  
 But in the latter case 
 $E_S$ would be an extension of ${\mathcal O}_S$ by ${\mathcal O}_S(1)$, 
 so this extension would split, thus contradicting the stability. 

 Therefore the Koszul complex of $s$ determines $E_S$ as an extension
 $$
 0 \rightarrow {\mathcal O}_S \rightarrow E_S 
 \rightarrow {\mathcal I}_Z(1) \rightarrow 0, 
 $$
 where ${\mathcal I}_Z$ is the ideal sheaf of $Z$, a zero 
 dimensional scheme of length six. Since $E_S$ is locally free, 
 Theorem 3.13 in \cite{L} implies that $Z$ has the following property:
 every hyperplane containing a colength one subscheme of $Z$
 contains $Z$. This means that $Z$ is of type $Z^k_6$, i.e. 
 $deg(Z) = 6$, $h^0(I_Z(1)) = 3+k$ and $h^1(I_Z(1)) = k$  for some $k\geq 1$. 

 \begin{lemma}\label{no-z26}
 Let $S$ be a smooth K3 surface in $G(2,6)$, with $Pic(S)=\ZZ\cO_S(1)$. 
Then
%\begin{enumerate} 
%\item $S$ contains no subscheme of type $Z_3^1$ or $Z_4^1$; 
%\item $S$ contains no subscheme of type $Z_6^k$ for $k>1$.
%\end{enumerate}

{(i)} $S$ contains no subscheme of type $Z_3^1$ or $Z_4^1$;

{(ii)} $S$ contains no subscheme of type $Z_6^k$ for $k>1$.
\end{lemma}

\proof {(i)} Since $S$ is an intersection of quadrics, 
if it contains a subscheme of type $Z^1_3$ then it contains the line 
spanned by this finite scheme. 
Next, suppose that $S$ contains a subscheme $Z$ of type $Z^1_4$.
Then $Z$ is the complete intersection 
of two quadrics in $\PP^2_Z = Span(Z)$,
and $\PP^2_Z \cap S = Z$ (otherwise $S$ would contain 
a conic or a line). 
But since $\PP^2_Z \cap G(2,6) = \PP^2_Z \cap S = Z$,
then $\PP^2_Z$ will be a purely 4-secant plane to $G(2,6)$
which is impossible -- any 4-secant plane to $G(2,6)$ 
must either lie in $G(2,6)$ or intersect $G(2,6)$ along  
a $1$-cycle of degree $2$. 

\smallskip
{(ii)} Since on $S$ there are no subschemes of type $Z^1_4$, then $S$ 
contains  
no subschemes of type $Z^2_5$ or $Z^3_6$. 
Next, suppose that on $S$ there exists a 
$0$-scheme $Z$ of type $Z^2_6$. Then $\PP^3_Z = Span(Z)$
will be a $3$-space in $\PP^8_S = Span(S)$ such that 
$\PP^3_Z \cap S \supset Z$. 

Consider a subscheme $Y\subset Z$ of length $5$. 
It must be of type $Z^1_5$. Consider the sheaf $E_Y$ 
obtained from $Y$ by the Serre construction, that is by 
the non trivial extension 
$$0\ra \cO_S\ra E_Y\ra I_Y(1)\ra 0. $$
This extension is unique since $Ext^1(I_Y(1),\cO_S)=H^1(I_Y(1))=\CC$.
Since $Y$ contains no subschemes of type $Z_4^1$, $E_Y$ is a vector
bundle \cite{L}. We check that it must be stable. If 
$L$ is a line bundle on $S$ with a non trivial sheaf homomorphism
$L\ra E_Y$, either the composition $L\ra I_Y(1)$ is non trivial 
or $L$ maps to $\cO_S\subset E_Y$; Since 
$Pic(S)=\ZZ\cO_S(1)$, in both cases we conclude 
that $L=\cO_S(\ell)$ for some $\ell\le 0$, which implies 
the stability. Using the argument of \cite{Mu4}, Theorem 4.5 (ii), we 
conclude that $E_Y$ is isomorphic to the restriction $T^*_S$
of the tautological vector bundle on the ambient Grassmannian. 
In particular $E_Y$ is generated by global sections. But 
then $I_Y(1)$ is also generated by global sections, which 
means that $S\cap\PP^3_Z=Y$, a contradiction. \qed 

\medskip

 Applying this result to our punctual scheme $Z$, which is 
 contained in a general K3 surface in $G(2,6)$, we conclude 
 that it must be of type $Z_6^1$, that is
 $h^0(I_Z(1)) = 4$, $h^1(I_Z(1)) = 1$ and $Z$ spans 
 a $\PP^4$ in $\PP^8_S = Span(S)$. 
 From the exact sequence 
$$0\ra H^1(E_S)\ra H^1(I_Z(1))\ra H^2(\cO_S)\ra 0$$ 
we deduce that $h^1(E_S) = 0$, hence $h^0(E_S)=5$.

Now (2) follows from the exact sequence 
 $$0 \rightarrow E(-1) \rightarrow E \rightarrow E_S \rightarrow 0,$$
 since by \ref{mod-5} we get $5\leq h^0(E)\le h^0(E_S)=5$. 
 This implies that $h^1(E)=0$, hence also $h^1(E(-1))=0$.
 \qed

 \begin{proposition}\label{glob-gen}
 Any vector bundle in $M_X(2;1,6)$ is globally generated.
 \end{proposition}

 \proof 
 Let $x\in X$ be any point. A crucial fact will be the following:

 \begin{lemma}\label{NL}
 In the linear system $|I_x(1)|$ of hyperplane sections of $X$ 
 passing through $x$, the set of smooth sections $S$ such that 
 $Pic(S)=\ZZ\cO_S(1)$ is a non empty open subset for the countable
 Zariski topology.
 \end{lemma}

\proof A Lefschetz pencil of
hyperplane sections of $X$ is defined 
by a line in the dual projective space,  cutting the dual variety
$X^*$ transversely at smooth points (see \cite{V}, Proposition 14.9).
Fix a point $x\in X$, and denote by $H_x$ the family of hyperplane
sections of $X$ containing $x$.  
The singular locus of $X^*$ cannot coincide with its intersection 
with the hyperplane $H_x$, so  there exists a Lefschetz pencil of hyperplane
sections of $X$ all passing through $x$. 
Then the proof of the
Noether-Lefschetz theorem by monodromy applies verbatim (see \cite{V},
Corollaire 15.28 and Th\'eor\`eme 15.33).\qed

\medskip
 Let $S \subset X$ be such a hyperplane section through $x$. 
 Since $h^1(E_S(-1))=0$ then $h^0(E)\simeq h^0(E_S)$. 
Thus to prove that $E$ is generated by global sections at $x$,
we just need to show that $E_S$ is generated by global sections. 
 Consider once again the Koszul complex 
 $$
 0 \rightarrow {\mathcal O}_S \rightarrow E_S 
 \rightarrow {\mathcal I}_Z(1) \rightarrow 0, 
 $$
 defined by some non trivial section of $E_S$. The bundle
 $E_S$ will be generated if and only if ${\mathcal I}_Z(1)$ is generated. 
 So we need to prove that $Z$ is cut out on $X$, scheme-theoretically,
 by its linear span -- in other words, that $Z$ cannot be contained 
 in a $Z'$ of type $Z_7^2$. 

 Suppose the contrary. 
 Then $h^1(I_{Z'}(1))=2$, and by 
 \cite{Mo}(page 22), $Z'$ defines a rank three sheaf $F$ over $S$
 as the universal  extension 
 $$
 0 \rightarrow O_S \otimes H^1(I_{Z'}(1)) \rightarrow F 
 \rightarrow I_{Z'}(1) \rightarrow 0.
 $$
 %
 %This sheaf is torsion free, and the induced cohomology sequence gives 
 
 Since by Lemma \ref{no-z26} $S$ does not contain any subscheme of type 
$Z^2_6$, then 
 any proper subscheme $Z'' \subset Z'$ has 
 $h^1(I_{Z''}(1)) \le 1$, and by \cite{Mo}(the Lemma on page 23) 
 the sheaf $F$ will be locally free, i.e. a vector bundle. 
 
 Next, consider the induced cohomology sequence 
 $$
 0\ra H^1(F)\ra H^1(I_{Z'}(1))\ra H^1(I_{Z'}(1))\ra H^2(F)\ra 0.
 $$
 Since $F$ has been defined as a universal extension, the
 middle map must be an isomorphism, hence $h^1(F)=h^2(F)=0$,  
 and $h^0(F)=6$; in particular $\chi(F) = 6$. 

 We shall see that $F$ is stable. Suppose it is not. 
 Since $c_1(F)$ is the hyperplane class of $S$ and $Pic(S)=\ZZ\cO_S(1)$, 
 then either there exists a non trivial map $\cO_S(1)\ra F$, 
 or there exists a non trivial map $F\ra\cO_S$. 
 But in the first case 
 the induced map to $I_{Z'}(1)$ vanishes and the map factors 
 through $\cO_S \otimes H^1(I_{Z'}(1))$, a contradiction. 
 Therefore there exists a non trivial map $F\ra\cO_S$, and denote 
 this map by $f$.  
 But then  the restriction of $f$ to $\cO_S \otimes H^1(I_{Z'}(1))$ 
 cannot be zero, since otherwise $f$ would descend to $I_{Z'}(1)$, 
 which is impossible. 
 %So we get a splitting of $O_S \otimes H^1(I_{Z'}(1))$
 %as $\cO_S\oplus\cO_S$, where the first factor maps to zero by $f$,
 %and the second isomorphically to $\cO_S$. In particular 
 Therefore $f$ is surjective, and since $h^2(\cO_S)=1$ then $h^2(F)>0$,
 again a contradiction. 

 %Thus $F$ is a stable torsion free sheaf of rank three
 Thus $F$ is a stable vector bundle of rank three
 with Chern numbers $c_1=1$ and $c_2=7$,  
 on the smooth $K3$ surface $S = S_{14}$ 
 with $Pic(S)=\cO_S(1)$.  
 In other words, the stable bundle $F$ belongs to the 
 Mukai's moduli space $M_S(r,L,s)=M_S(3,{\cO}_S(1),3)$
 of simple sheaves $E$ on $S$, with $rank(E) = r = 3$ 
 and $s = \chi(E) - r = 3$ (see above); 
 in particular this moduli space must be non-empty.   
 But by Theorem 0.1 in Mukai \cite{Mu1},
$M_S(3,{\cO}(1),3)$ 
 should have dimension = $deg\ S - 2.r.s = 14 - 2.3.3 + 2 = -2$,  
 contradiction. \qed

 \begin{corollary}\label{acm}
 Any vector bundle in $M_X(2;1,6)$ is arithmetically Cohen-Macaulay.
 \end{corollary}

 \proof We first check that $E(2)$ is Castelnuovo-Mumford regular,
 that is, $h^1(E(1))=h^2(E)=h^3(E(-1))=0$. By Serre duality, 
 $h^3(E(-1))=h^0(E^*)=h^0(E(-1))=0$. The vanishing of $h^2(E)$
 has already been checked in \ref{mod-5}, and
 $$H^1(E(1))=H^1(K_X\otimes E\otimes\det E\otimes\cO_X(1)),
 $$
 so the vanishing of $H^1(E(1))$ follows from 
 Griffiths' vanishing theorem, see e.g. \cite{Dem}.

 The Castelnuovo-Mumford regularity of $E(2)$ ensures that
 $$H^1(E(k+1))=H^2(E(k))=H^1(E(-k-2))=0 \quad \forall k\geq 0.
 $$
 Since $H^1(E)=H^1(E(-1))=0$ by Proposition \ref{mod-6}, 
 then $H^1(E(k))=0$ for all $k\in\ZZ$. Together with the 
 Serre duality, the last yields $H^2(E(k))=0$ for all $k\in\ZZ$.
 \qed

\medskip
Note that $E$ itself is not Castelnuovo-Mumford regular, 
even $E(1)$ is not since $h^3(E(-2))=h^0(E^*(1))=h^0(E)\ne 0$. 
That's what makes the proof of Proposition \ref{glob-gen} rather 
intricate.

 \subsection{Elliptic sextics}

 Since a vector bundle $E$ in $M_X(2;1,6)$ is globally 
 generated, a general section of $E$ will vanish along 
 a smooth curve $C$, an elliptic sextic in $G(2,6)$. 
 The associated Koszul complex is 
 $$0\ra\cO_X\ra E\ra I_C(1)\ra 0,$$
 and since $h^2(\cO_X)=0$  and $h^1(E)=0$ by Proposition \ref{mod-6},
 we deduce that $h^1(I_C(1))=0$. This means that $C$
 is {\it projectively normal}. 

 Consider the tautological rank two vector bundle $T^*$
 on $G(2,6)$, and restrict it to $C$. By the Atiyah classification 
 of rank two bundles on elliptic curves (see \cite{A}), and the fact 
 that $T^*$ is generated by global sections, 
 one of the following two possibilities must take place: 
 %we deduce that we must be in one of the following cases:
 
 {\it
 \smallskip
 \begin{enumerate}
 \item $T^*_{C}=L\op M$ is a direct sum of two line bundles
 of degrees $(3,3)$, $(2,4)$ or $(0,6)$;
 
 \smallskip
 \item $T^*_{C}=F\ot N$ is a direct product of a degree
 three line bundle $N$ with the unique vector bundle $F$ on $C$
 obtained as a non trivial self-extension of $\cO_C$.
 \end{enumerate}
 }
 
 \smallskip
 We say correspondingly that $C$ is split of type $(a,b)$, 
or unsplit.

If $C$ is split, in fact it cannot be of type $(0,6)$. Indeed,
this would imply that the restriction of $T$ to
 $C$ has a constant
 factor. This would mean that the lines in $\PP^5$ parameterized by 
 $C$ contain some fixed vector. But such lines are parameterized
 by a $\PP^4$ in $G(2,6)$, and we would conclude that $C$ is not 
 projectively normal, a contradiction. 

So the space of global sections
 $$H^0(C,T^*_{C})=H^0(C,L)\op H^0(C,M)$$
 has dimension $6$, which is the same as the dimension 
of $H^0(G(2,6),T^*)=V^*$.
 Note that the non-injectivity of the restriction map 
 $$H^0(G(2,6),T^*)\ra H^0(C,T^*_{C})$$ 
 implies that $C$ is contained 
 in a copy of $G(2,5)$. If this is not the case, $H^0(C,L)$ and 
 $H^0(C,M)$ can be identified with the orthogonal spaces 
 $A^{\perp}$ and $B^{\perp}$ to two supplementary spaces $A$ and $B$ 
 in $V$, and then the curve $C$ must be contained in 
 the intersection of $G(2,V)$ with the Segre variety 
 $$\Sigma_C=\PP A\times
 \PP B\subset\PP(A\ot B)\subset \PP(\wedge^2V).$$

 \begin{lemma}\label{sextype}
 Let $C$ be a smooth projectively normal elliptic sextic in $X$.
Then $C$  is split of type $(3,3)$, or unsplit.
 \end{lemma}

 \proof First we notice that $C$ cannot be contained in a 
copy $G(2,H)$ of $G(2,5)$. This is simply because $G(2,5)$
has degree 5, while $C$ has degree 6. So the intersection 
of $G(2,H)$ with $\PP^9_X$ cannot be proper. But then $X$
contains a surface of degree 5, contradicting the fact that 
its Picard group is generated by the hyperplane class -- or 
$X$ is contained in $G(2,H)$, which is clearly impossible 
for $X$ general. 

We must exclude the possibility that $C$ be of type $(2,4)$. 
Suppose it is. Then all the lines $\ell_x$, for $x\in C$, 
intersect a fixed line $D$ in $\PP^5$, identified with the 
image of $C$ by a complete linear system of type $g_2^1$. 
In particular, through a general point of $D$ pass two
lines $\ell_x$ and $\ell_{x'}$. But then 
the line $D$ must be contained in the singular locus $\Gamma(W)$
of the Palatani quartic (see 3.3), in contradiction with the fact that
$\Gamma(W)$ is a smooth curve of genus $26$. \qed 

 \medskip 
 A split elliptic sextic $C$ of type $(3,3)$ has two unisecant 
 planes, that is, two planes meeting all the lines parameterized
by $C$: with the notations above, they are the planes $\PP A$ 
and $\PP B$. In the unsplit case 
 there is a unique such plane. 
 
 We denote by $\cS(X)$ 
 the family of smooth projectively normal elliptic 
sextic curves $C \subset X$, an open subset of the Hilbert scheme of
 $X$. Let $\cS(X)_{un}$ be the subfamily of these $C \in \cS(X)$ that have  
 only one unisecant plane.

 \begin{lemma}\label{s1_5}
 For $X$ general, 
 the subfamily $\cS(X)_{un} \subset \cS(X)$ has dimension at most $5$.
 \end{lemma} 

 \proof 
 Let $C\subset G(2,6)$ be a smooth projectively normal elliptic
 sextic curve with only one unisecant plane. Then the restriction $T^*_C$
 of the tautological rank two bundle fits into an unsplit
 extension $0\ra L\ra T^*_C\ra L\ra 0$ for some degree three
 line bundle $L$ on $C$. This line bundle maps $C$ isomorphically
 to a plane cubic $E$ inside $\PP H^0(L)^*$. 

 Taking global sections, we get $0\ra H^0(L)\ra H^0(T^*_C)\ra 
 H^0(L)\ra 0$. If the restriction map $H^0(T^*)=V^*\ra 
 H^0(T^*_C)$ is not injective, then $C$ is contained in a $G(2,5)$. 
 Otherwise $H^0(T^*_C)=V^*$, and we deduce a projection map
 $\mu : V\lra H^0(L)^*\hookrightarrow V$. Let $\Lambda\subset V$
 be transverse to $M:=H^0(L)^*$, so that the restriction of $\mu$
 to $\Lambda$, that we denote by $\nu$, is an isomorphism.  
From the commutative diagram
$$\begin{array}{ccccc}
L^* & \ra & T_C & \ra & L^* \\
\downarrow & & \downarrow & & \downarrow \\
M\ot\cO_C & \ra & V\ot\cO_C & \ra & M\ot\cO_C, 
\end{array}$$
whose vertical maps are injective, 
we deduce that over a point $[e]\in E\simeq C$, where $e\in M$, 
the two-dimensional subspace $T_C$ of $V$ is generated
by $e$ and a vector $f$ mapping to $e$ by $\mu$. So $f$ must be 
of the form $\nu^{-1}(e)+\theta(e)$ for some $\theta(e)\in M$, defined up to 
translation by some multiple of $e$. Thus $\theta$ has to be interpreted
as a global section of $\mathcal{Hom}(L^*,M\ot\cO_C/L^*)$. 
From the exact sequence of vector bundles $0\ra\cO_C\ra M\ot L\ra 
\mathcal{Hom}(L^*,M\ot\cO_C/L^*)\ra 0$ we deduce the sequence 
$$0\ra\CC\ra End(M)\ra Hom(L^*,M\ot\cO_C/L^*)\ra\CC\ra 0.$$  
In particular $Hom(L^*,M\ot\cO_C/L^*)$ has dimension 9. 

 We can now count the number of parameters for $C$: we have 
 9 parameters for the three-space $H^0(L)^*\subset V$, 9
 parameters for the cubic curve $E\subset \PP H^0(L)^*$, 
 then 9 again for $\mu$ and 9 for $\theta$, minus one since 
 multiplying $\mu^{-1}$ and $\theta$ by a same scalar does not
 change the resulting curve: this makes 35 parameters. Each 
 of our curves $C$ spans a $\PP^5$, and there is a 20-dimensional
 family of $\PP^9$'s containing this $\PP^5$. Since the Grassmannian
 of $\PP^9$'s in $\PP^{14}$ has dimension 50, our claim follows. 
 \qed 

 \begin{corollary}\label{sexdef}
 Any smooth projectively normal elliptic sextic in $X$
 with only one unisecant plane can be deformed into an 
 elliptic sextic with two unisecant planes.
 \end{corollary}

 \proof Let $C\in\cS(X)_{un}$. The dimension at
 $C$ of the family $\cS(X)$ is at least 
 $\chi(N_{C/X})$ (see e.g. \cite{Se}, Corollary 8.5), 
which  by Riemann-Roch is equal to
 $deg(C)=6$. So by Lemma \ref{s1_5}, the curve $C$ can be 
 deformed outside $\cS(X)_{un}$.\qed

\medskip\noindent {\bf Definition}. 
{\it Let $M_X(2;1,6)^0$ denote the
open subset of vector bundles in $M_X(2;1,6)$ having 
a section whose zero locus is a split elliptic sextic.}

\medskip By the previous statement any vector bundle in 
$M_X(2;1,6)$ is in the closure of the open set $M_X(2;1,6)^0$. In fact
we will finally conclude that any vector bundle in 
$M_X(2;1,6)$ belongs to $M_X(2;1,6)^0$, see Theorem \ref{main}. 

 \medskip
 Let $C\in\cS(X)$ be a split elliptic sextic, with 
 two unisecant planes $\PP A$ and $\PP B$. The linear span 
 of the Segre variety $\Sigma_C=\PP A\times \PP B$ is 
 $\PP(A\ot B)\subset \PP(\wedge^2V)$, whose orthogonal 
in the dual space is 
 $$\PP(\wedge^2A^{\perp}\oplus \wedge^2B^{\perp})
 \subset \PP(\wedge^2V^*).$$
 Since the span of $X$ meets $\PP(A\ot B)$ along the five dimensional
linear span of $C$, 
 their orthogonal spaces meet each other along a line 
 $\ell\subset Y=X^{\perp}$. Note that $\ell$ only depends on the 
vector bundle $E$ defined by $C$, since the zero-loci of the sections 
of $E$ are parameterized by the projective space of its global sections, 
while $F(Y)$, being embedded in the abelian variety $F(Y)$,  contains 
no rational curve. We thus get a map $M_X(2;1,6)^0\rightarrow F(Y)$. 
An inverse of this map
will be constructed in the next section. 

\smallskip
We now use what we learned on elliptic sextics to prove that 
$M_X(2;1,6)^0$ is a smooth open subset of the moduli space. 

 \begin{proposition}\label{vb-smooth}
 Let $E\in M_X(2;1,6)^0$ be a vector bundle. Then
 $$h^0(End(E))=1,\quad h^1(End(E))=2,\quad h^2(End(E))=h^3(End(E))=0.$$
 In particular $M_X(2;1,6)$ is smooth at $E$.
 \end{proposition}

 \proof By Serre duality $h^2(End(E))=h^1(End(E)(-1))=h^1(E\ot E(-2))$.
 Consider the Koszul complex
 $$0\ra\cO_X\ra E\ra I_C(1)\ra 0$$
 of a section of $E$ vanishing along a smooth elliptic sextic $C$
 with two unisecant planes. Twisting this sequence by $E(-2)$ and using 
 \ref{acm}, we are reduced to prove that $h^1(E\ot I_C(-1))=0$, thus 
 that $h^0(E_C(-1))=h^0(E_C^*)=0$. 

 Suppose the contrary. Since $E$, hence $E_C$, is globally generated,
 we deduce that $E_C$ splits. Since it is isomorphic to the normal
 bundle of $C$ in $X$, we get that  $N_{C/X}=\cO_C\op\cO_C(1)$. Now
 $C$ is a linear section of a Segre variety $\Sigma_C=\Sigma\simeq
\PP A\times \PP B$ by the linear span of $X$, meeting the linear span of 
 $\Sigma$ along the span of $C$. Let $G$ denote the
 intersection of the Grassmannian with the linear span of $X$ and
 $\Sigma$, a codimension two linear space.
 We have the following diagram of normal bundles on $C$:
 $$\begin{array}{ccccc}
 & & 0 & & \\ & & \uparrow & & \\
  & & N_{X/G}|_C = \cO_C(1)^{\op 3} & & \\
  & & \uparrow & & \\
0\ra N_{C/\Sigma}=\cO_C(1)^{\op 3} & \ra & N_{C/G} & \ra & 
N_{\Sigma/G}|_C\ra 0 \\
  & & \uparrow & & \\
  & & N_{C/X} & & \\
  & & \uparrow & & \\
  & & 0& & 
 \end{array}$$
 The induced map from $N_{C/\Sigma}$ to $N_{X/G}|_C$ is an isomorphism,
 hence also the induced map from $N_{C/X}$ to $N_{\Sigma/G}|_C$, 
 and we get an injective map from  $N_{C/X}$ to $N_{\Sigma/G(2,6)}$   
 restricted to $C$. 
 Remember that $T^*_C=L\op M$ for some line bundles $L$
 and $M$ of degree three, whose spaces of global sections 
are $H^0(L)=A^{\perp}$ and $H^0(M)=B^{\perp}$. In particular 
$\cO_C(1)=L\ot M$. We compute 
 that the normal bundle of $\Sigma$ inside $G(2,6)$ is
 $$N_{\Sigma/G(2,6)}=(L\ot B/L^*)\oplus (M\ot A/M^*).$$
 Suppose we have a non trivial morphism from $\cO_C(1)$, which we
 supposed to be a factor of $N_{C/X}$, to one of these direct factors,
 say the first one. Then we deduce that $Hom(M,B/L^*)\neq 0$. But the
 rank two vector bundle $B/L^*$ on $C$ is globally generated and has
 degree three, and from the Atiyah classification we deduce that 
 $B/L^*=M\op\cO_C$. Then the restriction to $C$ of the
 tautological quotient bundle $Q_C$, which is isomorphic to
 $B/L^*\op A/M^*$, has a trivial factor. This implies that the 
 linear span of $C$ is contained in a copy of $G(2,5)$, a contradiction.
 
That the vanishing of $h^2(End(E))$ implies the smoothness of 
$M_X(2;1,6)$ at $E$ is well-known, see e.g. \cite{HL}, Theorem 4.5.4. 
\qed

 \subsection{Non locally free sheaves in $M_X(2;1,6)$}

In this section we give a complete description of the non 
locally free sheaves in our moduli space $M_X(2;1,6)$. A 
similar study has been made by Druel on the cubic threefold. 
Our discussion follows closely that of \cite{D}, but 
Mukai's theorems on Fano threefolds and K3 surfaces of genus 8
will play a crucial role. 

 \begin{proposition}\label{mod-12}
 If $E \in M_X(2;1,6)$ is not locally free, 
 there exists a unique line $\ell \subset X$ such that $E$ 
 fits into an exact sequence
 $$0\ra E\ra T^*_X \rightarrow \cO_{\ell} \rightarrow 0.$$ 
 \end{proposition}

 Note that $T^*$ restricted to $\ell$ is isomorphic to 
 $\cO_{\ell}\op\cO_{\ell}(1)$, so that the rightmost arrow
 is uniquely defined up to scalar. Thus the line $\ell$ defines
 $E$ uniquely. 
%In particular, we get a map $M_X(2;1,6)^1\rightarrow F(X)$.

 \proof
 Let $F = E^{**}$ denote the bidual of $E$, and let $R=F/E$.
 Since $F$ is reflexive, the singular locus $S(F)$ has codimension
 at least three, see Lemma 1.1.10 in Ch.2 of \cite{OSS}. 
 Therefore the restriction $F_S$ of $F$ to the
 general hyperplane section $S$ of $X$ is locally free, and
 the restriction $R_S$ of $R$ has finite support.  
 We have $c_1(F_S)= 1$ and  $c_2(F_S) = 6 - \mathrm{length}(R_S)$,
 so by Riemann-Roch 
 $\chi(F_S)=2\chi(\cO_S)+(c_1(F_S)^2-2c_2(F_S))/2 =\mathrm{length}(R_S)+5$. 

 This rank two vector bundle $F_S$ is semistable by Maruyama's
 restriction theorem. By Mukai (\cite{Mu1}, Theorem 0.1), the
 moduli space of simple sheaves on $S$ with these invariants 
 is smooth at $F_S$ and has dimension $c_1(F_S)^2-4\chi(F_S)+10=
 4-4\mathrm{length}(R_S)$. In particular 
 $$\mathrm{length}(R_S)\leq 1.$$ 

 \smallskip
First, we shall see that $\mathrm{length}(R_S)$ cannot be zero.  
 Suppose the contrary. 
 Then $R$ has finite support and $c_2(F)=6$. 
 Moreover, by Riemann-Roch $\mathrm{length}(R)=\chi(F(-1))=c_3(F)/2$. 
 Since $\chi(F_S)=5$ and $h^2(F_S)=h^0(F_S^*)=0$ by stability,
 the vector bundle $F_S$ has a non trivial section. The corresponding
 Koszul complex gives 
 $$0\ra \cO_S\ra F_S\ra I_Z(1)\ra 0,$$
 where the finite scheme $Z$ has type $Z_6^1$. The associated
 long exact sequence 
 $$0\ra H^1(F_S)\ra H^1(I_Z(1))\ra H^2(\cO_S)\ra 0$$
 gives $H^1(F_S)=0$. Moreover, $H^1(F_S(k))=H^1(I_Z(k+1))$ 
 for $k>0$. But this is zero, since the fact that $Z$ has type $Z_6^1$ 
 easily implies that the restriction map $H^0(\cO_S(k+1))
 \ra H^0(\cO_Z(k+1))$ is surjective. Now the exact sequence 
 $$0\ra F(-1)\ra F\ra F_S\ra 0$$
twisted by $\cO_X(k)$,
 together with the vanishing of $H^1(F_S(k))$ for all $k\geq 0$,
 imply that $H^2(F(k-1))$ embeds inside $H^2(F(k))$,
 hence it is always zero since it certainly vanishes for $k$ large enough.
 In particular $h^2(F(-1))=0$, and since $h^0(F(-1))=0$ by stability,  
 we conclude that $\chi(F(-1))\leq 0$. Hence $c_3(F)=0$, so that 
 $F$ is locally free by \cite{H2},   and $\mathrm{length}(R)=0$; 
so that $E$ is isomorphic to  $F$. But then $E$ is 
locally free -- contradiction.

 \smallskip
 Therefore the only possibility left is $\mathrm{length}(R_S) = 1$, so 
$R_S=\cO_p$ for some point $p\in S$. 
 In this case, $R$ must be supported on  the union 
 of a line $\ell$ with a finite set, and must have multiplicity 
 one on that line.
 By Riemann-Roch $\chi(R(-1))=\chi(F(-1))
 =c_3(F)/2$. The restriction $F_S$ has Chern classes $c_1(F_S)=h$
 and $c_2(F_S)=5$. Therefore by \cite{Mu4}, $F_S$ is the same as 
 the restriction $T_S^*$ to $S \subset G(2,6)$ of the tautological 
 bundle $T^*$ on $G(2,6)$.  
 Since $S$ is a linear section of $G(2,6)$, 
 we have a Koszul complex
 $$0\ra T^*(-6)\ra\cdots\ra T^*(-1)^{\op 6}\ra T^*\ra F_S\ra 0.$$
 Thus $h^1(F_S(k))=0$ as soon as $h^{q+1}(T^*(k-q))=0$ for any
 $q=0,\ldots ,6$. By Bott's theorem this holds true for any 
 $k\in\ZZ$. As in the previous case we conclude that $h^2(F(-1))=0$,
 hence $\chi(F(-1))\leq 0$, hence $c_3(F)=0$ and $F$ is locally free.
 Since $c_1(F)=h$ and $c_2(F)=5$, then again by \cite{Mu4} 
(see e.g. Theorem 1.10), 
 the bundle $F$ must be the same as $T^*_X$, the restriction to 
 $X$ of the tautological bundle $T^*$ on $G(2,6)$. 
 %Using Bott's theorem as before we deduce that $h^1(F(-1))=0$. 

 On the other hand, from the exact sequence 
 $$0\ra E_S\ra T_S^*\ra\cO_p\ra 0,$$ 
 we deduce that $h^1(E_S(k))=0$
 for $k\ge 0$, hence $h^2(E(-1))=0$, and then $h^1(E(-1))=0$
 since $\chi(E(-1))=0$, $h^0(E(-1))=0$ by stability and 
 $h^3(E(-1))=h^0(E^*)=0$ by Serre duality for sheaves. Thus
 $h^0(R(-1))=0$. Since Riemann-Roch gives
 $\chi(R(n))=n+1$, we conclude by \cite{D}, Lemme 3.2, 
that $R=\cO_{\ell}$. \qed

 \begin{proposition}\label{m1-smooth}
 Let $E \in M_X(2;1,6)$ be a non locally free sheaf. Then
 $$hom(E,E)=1, \quad ext^1(E,E)=2, \quad ext^2(E,E)=ext^3(E,E)=0.$$
In particular the moduli space $M_X(2;1,6)$ is smooth at $E$.
 \end{proposition}

 \proof 
 Since $E$ is stable, $hom(E,E)=1$. By Serre duality we deduce 
 that $ext^3(E,E)=hom(E,E(-1))$
 must vanish, since otherwise we would clearly get $hom(E,E)>1$.
 By Riemann-Roch $hom(E,E)-ext^1(E,E)+ext^2(E,E)-ext^3(E,E)=-1$,
 so we just need to check that $ext^1(E,E)=2$. 

 By \ref{mod-12} there is a line $\ell$ on $X$ such that $E$ fits 
 into an exact sequence 
 $$0\ra E\ra T^*_X \rightarrow \cO_{\ell} \rightarrow 0.$$
We have $ext^i(T_X^*,T_X^*)=0$ for $i>0$: use the Koszul complex
of $X$ and Bott's theorem as in 3.1 (this means that the moduli 
space $M_X(2;1,5)$, which by Mukai's theorem reduces to one point, 
 is smooth).  
Tensoring the previous short exact sequence by 
 $T_X$ and taking cohomology, we deduce that $h^i(T_X\ot E)=0$
 for $i>0$. 
 Applying the functor $Hom(.,E)$ to the same short sequence, 
 we deduce that 
 $$Ext^1(E,E)\simeq Ext^2(\cO_{\ell},E).$$
 Now we apply $Hom(\cO_{\ell},.)$ and obtain
 $$\begin{array}{l}
 Ext^1(\cO_{\ell},T_X^*)\ra Ext^1(\cO_{\ell},\cO_{\ell})\ra \\
 \hspace*{3cm}
 \ra Ext^2(\cO_{\ell},E)\ra Ext^2(\cO_{\ell},T_X^*)
 \ra Ext^2(\cO_{\ell},\cO_{\ell})
 \end{array}$$
 To compute $Ext^i(\cO_{\ell},T_X^*)$, recall that is can be obtained
 as the abutment of the spectral sequence with order two terms
 $$E_2^{j,i-j}=H^j({\mathcal Ext}^{i-j}(\cO_{\ell},T_X^*))
 =H^j({\mathcal Ext}^{i-j}(\cO_{\ell},\omega_X)\ot T_X^*(1))).
 $$
 Since $\ell$ is smooth of codimension two, the sheaf 
 ${\mathcal Ext}^k(\cO_{\ell},\omega_X)=0$ for $k<2$, and 
 ${\mathcal Ext}^2(\cO_{\ell},\omega_X)=\omega_{\ell}$.
 We deduce that $Ext^1(\cO_{\ell},T_X^*)=0$ and 
 $$Ext^2(\cO_{\ell},T_X^*)=H^0({\mathcal Ext}^2(\cO_{\ell},\omega_X)
 \ot T^*_X(1)) =H^0(T^*_{\ell}(-1))=\CC.$$
 Now recall that $\Gamma(X)$, the family of lines on $X$ is a smooth 
 irreducible curve, see \ref{xlines}, and that 
 $ext^1(\cO_{\ell},\cO_{\ell})=1$ and $ext^2(\cO_{\ell},\cO_{\ell})=0$
for any line $\ell$. 
 Putting all this in the long exact sequence above, we finally 
 get that 
 $$ext^1(E,E)=ext^2(\cO_{\ell},E)=ext^1(\cO_{\ell},\cO_{\ell})
 +ext^2(\cO_{\ell},T_X^*)=1+1=2,$$
 which concludes the proof. \qed 

 \section{Vector bundles and projections}

Every vector bundle in $M_X(2;1,6)$ can be obtained by the 
Serre construction from a smooth elliptic sextic in $X$. In this
section we give an alternative construction. 

 \subsection{Projections associated to lines}

 Let $e_0,\ldots ,e_5$ be a basis of $V=\CC^6$, and let $f_0,
 \ldots ,f_5$ be the dual basis. Consider the A-line 
$\ell$ in $\PP(\Lambda^2V^*)$ generated by 
 \begin{eqnarray} \nonumber
  & f_0\we f_2+f_1\we f_3, \\ \nonumber
  & f_0\we f_4+f_1\we f_5.  
 \end{eqnarray}
 Its orthogonal $\ell^{\perp}$ is a special codimension two subspace in 
 $\PP^{14}=\PP(\Lambda^2V)$. 

 \begin{lemma}
 The singular locus of the linear section $G(2,6)\cap \ell^{\perp}$ is 
 a smooth plane conic $q^{\ell}$,
 parameterizing the singular points of the sections of $G(2,6)$ by 
 a hyperplane in the pencil $\ell$. The projective span of $q^{\ell}$ is 
 a plane $\pi^{\ell}$ which is not contained in $G(2,6)$.  
 \end{lemma}

 \proof This singular locus is the set of points $x\we y$ on $G(2,6)\cap
 \ell^{\perp}$
 at which the tangent space to $G(2,6)$ is not transverse to
 $\ell^{\perp}$. 
 This means that 
 one of the linear forms $f_0\we f+f_1\we f'$ of $\ell$ vanishes on the 
 affine tangent space, 
 which is $x\we V+y\we V$. This is possible only if the plane $\langle
 x,y
 \rangle$ is orthogonal
 to the linear forms $f_0,f_1,f,f'$. Note that this defines uniquely the 
 corresponding 
 point of $G(2,6)$, as the singular point of the intersection of $G(2,6)$ 
 with the hyperplane
 orthogonal to $f_0\we f+f_1\we f'$. Explicitly, letting $f=sf_2+tf_4$ 
 and $f'=sf_3+tf_5$, 
 this singular point is given by $x=te_2-se_4$ and $y=te_3-se_5$, so 
 that $x\we y=
 t^2e_2\we e_3-st(e_2\we e_5-e_3\we e_4)+s^2e_4\we e_5$ describes a 
 smooth conic, as claimed. 
 \qed

 \medskip
 Now we project 
 $G(2,6)\cap \ell^{\perp}\subset \ell^{\perp}$ 
 linearly from the projective plane $\pi^{\ell}$. The image 
 of this projection is a six-dimensional
 variety $G_{\ell}\subset\PP^9$. 

 \begin{proposition}
 The variety $G_{\ell}\subset\PP^9$ is projectively equivalent to the 
 Grassmannian $G(2,5)$ in the Pl\"ucker embedding.
 \end{proposition}

 \proof 
Keeping the previous notations, 
 the projective plane $\pi^{\ell}$ is generated by $e_2\we e_3, e_2
 \we e_5-e_3\we e_4$ and $e_4\we e_5$. 
 To describe the variety $G_{\ell}$, we first parameterize an open subset 
 of $G(2,6)$ as follows: 
 any plane transverse to $\langle e_0,e_1,e_3,e_5\rangle$ has a unique 
 basis of the form
 \begin{eqnarray} \nonumber
  & e_2+\a_0e_0+\a_1e_1+\a_3e_3+\a_5e_5, \\ \nonumber
  & e_4+\b_0e_0+\b_1e_1+\b_3e_3+\b_5e_5. 
 \end{eqnarray}
 This gives affine coordinates on $G(2,6)$. Let
 $\d_{ij}=\a_i\b_j-\a_j\b_i$. 
 The linear
  section $G(2,6)\cap \ell^{\perp}$ is defined by the conditions
 $$\b_0-\d_{13}=\a_0+\d_{15}=0.$$
 The image $G_{\ell}$ of the projection of this section along $\pi^{\ell}$ 
 is, in suitable coordinates, 
 the closure of the set of points in $\PP^9$ of the form
 $$ [1,\a_1,\b_1,\a_3+\b_5,\d_{01}, \d_{03}, \d_{05}, \d_{13}, \d_{15}, 
 \d_{35}].$$
 We can already say that this is a unirational variety, locally
 parameterized
  by $\a_1$,
 $\a_3$, $\a_5$, $\b_1$, $\b_3$, $\b_5$. Indeed, $\a_0$ and $\b_0$ are functions of
 these 
 parameters, 
 hence also $\d_{01}, \d_{03}, \d_{05}$. 

 Let us denote by $X,Y,Z,T,U_{01},U_{03}, U_{05}, U_{13}, U_{15}, U_{35}$ 
 our homogeneous 
 coordinates on $\PP^9$. 

 \begin{lemma} 
 $G_{\ell}$ has five quadratic equations, explicitly given by
 \begin{eqnarray}
 XU_{01}+YU_{13}+ZU_{15}=0, \\ 
 XU_{03}+ZU_{35}+TU_{13}=0, \\ 
 XU_{05}+YU_{35}-TU_{15}=0, \\ 
 YU_{03}+ZU_{05}-TU_{01}=0, \\ 
 U_{01}U_{35}-U_{01}U_{15}+U_{05}U_{13}=0.  
 \end{eqnarray}
 \end{lemma}

 \noindent {\it Proof of the lemma}. 
 We first note that $\d_{01}=\a_0\b_1-\b_0\a_1=-\d_{15}\b_1-\d_{13}\a_1$, 
 which gives $(1)$.
 We also have 
 \begin{eqnarray} \nonumber
 \d_{03}+\d_{15}\b_3+\d_{13}\a_3=0, \\ \nonumber
 \d_{05}+\d_{15}\b_5+\d_{13}\a_5=0.
 \end{eqnarray}
 Combining these identities to the Pl\"ucker relations 
 \begin{eqnarray} \nonumber
 \d_{13}\b_5-\d_{15}\b_3+\d_{35}\b_1=0, \\ \nonumber
 \d_{13}\a_5-\d_{15}\a_3+\d_{35}\a_1=0, 
 \end{eqnarray}
 we get the two equations 
 \begin{eqnarray} \nonumber
 \d_{03}+(\a_3+\b_5)\d_{13}+\b_1\d_{35}=0, \\ \nonumber
 \d_{05}+(\a_3+\b_5)\d_{15}+\a_1\d_{35}=0,
 \end{eqnarray}
 which are $(2)$ and $(3)$. Two other  Pl\"ucker relations are 
 \begin{eqnarray} \nonumber
 \d_{15}\b_0+\d_{01}\b_5-\d_{05}\b_1=0, \\ \nonumber
 \d_{01}\a_3-\d_{03}\a_1+\d_{13}\a_0=0. 
 \end{eqnarray}
 Adding them, we get $(4)$. Finally $(5)$ is itself a Pl\"ucker relation. \qed

 \medskip To conclude the proof of the Proposition, we need to check that these
 quadrics are the Pl\"ucker equations of $G(2,5)$ in a slightly disguised form. 
 This will imply that the 6-dimensional variety $G_{\ell}$ is contained
 in, 
 hence equal to,
 a copy of $G(2,5)$. 

 \smallskip We make the following substitution:
 $$\begin{array}{lll}
 \d_{01}=\Delta_{12}\quad &  \d_{13}=\Delta_{14}\quad &  X=\Delta_{45} \\  
 \d_{03}=\Delta_{13} &  \d_{15}=\Delta_{24} &  Y=-\Delta_{25} \\  
 \d_{05}=\Delta_{23} &  \d_{35}=\Delta_{34} &  Z=\Delta_{15} \\  
   &  &  T=-\Delta_{35}. 
 \end{array}$$
 The five quadrics of the Lemma become
 \begin{eqnarray} \nonumber
 \Delta_{12}\Delta_{45}-\Delta_{14}\Delta_{25}+\Delta_{15}\Delta_{24}=0,
  \\ 
 \nonumber 
 \Delta_{13}\Delta_{45}+\Delta_{15}\Delta_{34}-\Delta_{14}\Delta_{35}=0,
 \\ 
 \nonumber 
 \Delta_{23}\Delta_{45}+\Delta_{25}\Delta_{34}-\Delta_{24}\Delta_{35}=0,
 \\ 
 \nonumber 
 \Delta_{13}\Delta_{25}+\Delta_{15}\Delta_{23}-\Delta_{12}\Delta_{35}=0,
 \\ 
 \nonumber 
 \Delta_{12}\Delta_{34}-\Delta_{13}\Delta_{24}+\Delta_{14}\Delta_{23}=0. 
 \end{eqnarray}
 The proof is complete. \qed

 \subsection{Induced vector bundles}

 Let as above $X$ be a general prime Fano threefold of index one and 
 of genus $8$, and let $Y$ be its orthogonal cubic threefold. 
 For an A-line $\ell \subset Y$, 
 consider the linear projection 
 $$f_{\ell} : 
 G(2,6)\cap \ell^{\perp}\dashrightarrow G_{\ell}\simeq G(2,5).$$
 It can be easily seen that $f_{\ell}$ is a birational map,
 which is regular outside the conic $q^{\ell}$ swept out 
 by the singular points of the hyperplane sections 
 $H_y \subset G(2,6)$ defined by the points $y \in \ell$.

 The Fano threefold $X$ in $G(2,6)$ can't meet the conic
 $q^{\ell}$, otherwise $X$ would be singular at its 
 intersection points with $q^{\ell}$. 
 Therefore the restriction  
 $$f_{\ell,X} : X\lra G_{\ell}$$
 of $f_{\ell}$ to $X$ is regular (and a birational map to its image). 
 In particular, the pull-back of the tautological bundle 
 $T^*$ on $G(2,5) = G_{\ell}$ restricts to a rank two vector 
 bundle $E_{\ell}$ on $X$. 

\begin{lemma}\label{chern}
We have  $c_1(E_{\ell})=h$ and $c_2(E_{\ell})=6$.
\end{lemma}

\proof 
Clearly $c_1(E_{\ell})=h$ since $c_1(T^*)$ is the hyperplane class
of $G(2,5)$ and $E_{\ell}$ is its pull-back by a linear projection. 

Let $A,B$ denote complementary three-spaces in $V$ such that 
the line $\ell\subset\PP(\we^2A^{\perp}\op \we^2B^{\perp})$. 
Then $\ell^{\perp}$
 contains the Segre variety $\Sigma=\PP A\times \PP B$, which meets 
 the plane $\pi^{\ell}$ precisely along the conic $q^{\ell}$. 
 The image of $\Sigma$ by the projection $f_{\ell}$ is thus a 
 four dimensional quadric $Q\subset G_{\ell}$. Such a quadric is 
 a copy of $G(2,4)$ inside $G_{\ell}\simeq G(2,5)$. In particular,
 it can be described as the zero locus of a general section of the 
 tautological bundle $T^*$ on $G_{\ell}$.  
The induced section of $E_{\ell}$ vanishes along $\Sigma\cap X$, and 
this intersection is generically transverse; hence $\Sigma\cap X$ 
represents the second Chern class of $E_{\ell}$. 
Since $\Sigma$ has degree 6, we conclude that $c_2(E_{\ell})=6$. \qed

 \begin{lemma}
 The vector bundle $E_{\ell}$ is stable. 
 \end{lemma}

 \proof We just need to check that there is no embedding 
 $\cO_X(1)\hookrightarrow E_{\ell}$, or equivalently that 
 $h^0(E_{\ell}(-1))=h^0(E_{\ell}^*)=0$. But since $E_{\ell}$
 is globally generated, a non zero section of $E_{\ell}^*$
 cannot vanish, and we would deduce that $E_{\ell}=\cO_X
 \op\cO_X(1)$. Then the tautological bundle restricted
 to $f_{\ell}(X)$ would have a trivial factor, which would 
 mean that $f_{\ell}(X)$ is a $\PP^3$ inside 
 $G_{\ell}$ -- a contradiction. 
 \qed

 \medskip
 Now let $E$ in $M_X(2;1,6)^0$ be any vector bundle. A general
 section of $E$ vanishes along a smooth split elliptic sextic $C$ with
 two unisecant planes. These two planes define complementary
 three-spaces $A,B$ in $V$, and $\PP(\wedge^2A^{\perp}\oplus
\wedge^2B^{\perp})$
 cuts the cubic $Y$ orthogonal to $X$ along an A-line $\ell$. 

But by the proof of Lemma \ref{chern}, we know that the 
vector bundle $E_{\ell}$ has a section whose zero locus is
precisely $C$. Since there is a unique vector bundle in $M_X(2;1,6)$
 with a section vanishing along $C$, we get $E_{\ell}\simeq E$ 
 and we deduce:

 \begin{proposition}\label{bij}
 The correspondence $\ell\mapsto E_{\ell}$ defines a bijection 
 between A-lines  in the orthogonal cubic $Y$, and vector
 bundles in $M_X(2;1,6)^0$.
 \end{proposition}

\begin{corollary}\label{irred}
$M_X(2;1,6)$ is an irreducible surface. 
 \end{corollary}

\proof Non locally free sheaves in $M_X(2;1,6)$ describe a curve 
contained, by \ref{m1-smooth}, in the smooth locus of a two-dimensional 
component of $M_X(2;1,6)$. In particular they belong to the closure of
the open subset of vector bundles, hence by \ref{sexdef} to the 
closure of $M_X(2;1,6)^0$. Finally, Proposition \ref{bij} implies that 
$M_X(2;1,6)^0$ is irreducible, and by \ref{vb-smooth} it is a 
smooth surface. \qed

 \section{The main result}

 Let again $J(X)$ denote the intermediate Jacobian of $X$. As in 
 \cite{D}, Theorem 4.8,  once we fix a base point $E_0\in M_X(2;1,6)$,
the second Chern class defines 
 a morphism $$c_2 : M_X(2;1,6)\lra J(X),$$
which is 
 uniquely defined up to a translation. 

Consider a vector bundle $E_{\ell}\in M_X(2;1,6)^0$. A general 
section of $E_{\ell}$ vanishes along a smooth split elliptic sextic
$C$, and $c_2(E)=[C]$. Let $\PP A$ be one of the two unisecant
planes to $C$ -- we can suppose that the line $\ell$ can be defined
as at the beginning of 6.1, where the dual basis $e_0,\ldots ,e_5$
of $f_0,\ldots ,f_5$ is such that $A=\langle e_0,e_2,e_4\rangle$.
Consider the intersection with $X$
of the Schubert cycle $\sigma_{20}(\PP A)$, the degree 
9 cycle of lines meeting $\PP A$. This intersection is the 
union of $C$ with  a cubic curve $D$. A
point  in $D$ but not in $C$ is defined by a tensor $\omega
=\alpha+\beta$, where $\alpha\in\wedge^2A$ and $\beta\in A\ot B$. 
Since $X\subset\ell^{\perp}$, $\alpha$ must be a non zero 
multiple of $e_2\wedge e_4$. Then for $\omega$ to have rank two, 
$\beta$ must not involve $e_0$. In particular, $D$ is contained
in $G(2,H)\simeq G(2,5)$, if  $H$ is the hyperplane generated
by $e_1,\ldots ,e_5$. More precisely, we can write $\omega=
e_2\wedge e_4+e_2\we u_2+e_4\we u_4$, with $u_2,u_4\in B$. But 
then $u_2$ and $u_4$ must be collinear, which implies that the 
line $\langle e_2,e_4\rangle$ in $\PP^5$ meets the line defined 
by $\omega$. We deduce:

\begin{lemma}
Let $C$ be a projectively normal elliptic sextic in $X$, 
with a unisecant plane $\PP A$. Let $D$ be the residual 
cubic curve of $C$ in the Schubert cycle $\sigma_{20}(\PP A)\cap X$.
 
Then there is a line $ax(D)$ and a hyperplane $H$ in $\PP^5$, 
with $ax(D)\subset H$, such that $D$ is contained in $G(2,H)$
and $ax(D)$ is a unisecant line of $D$.
\end{lemma}

Of course we have checked this only for a split elliptic 
sextic, but by continuity the statement holds also in the
unsplit case. We call the line $ax(D)$ the {\it axis} of the 
cubic curve $D$. 

\medskip 
Now we notice that $G(2,5)$ has degree $5$, and contains
both the cubic $D$ and the conic  $q_{\ell}$. Thus the intersection
$G(2,H)\cap X=D\cup q_{\ell}$ and, up to some fixed translations, we
have
$$c_2(E_{\ell})=[C]=-[D]=[q_{\ell}]\in J(X).$$

 \begin{theorem}\label{main}
 The map $c_2$ defines an isomorphism between the moduli
 space $M_X(2;1,6)$ and a translate in $J(X)$ of the Fano
 surface $F(X)$ of conics in $X$.
 \end{theorem}

 \proof We have just seen that the map $c_2$ sends $M_X(2;1,6)^0$
to (a translate of) the open set of conics $q_{\ell}$ 
in $F(X)$, which by \ref{-1} corresponds to the set of  A-lines in $F(Y)$;
and by \ref{bij} this is a bijection. Since $M_X(2;1,6)$ is 
irreducible by \ref{irred}, we deduce that $c_2$ maps 
$M_X(2;1,6)$ to $F(X)$ in $J(X)$. 
Moreover, it follows from \ref{mod-12} that non locally free sheaves
in $M_X(2;1,6)$ are mapped bijectively to the image of 
the curve $\Gamma(X)$ of lines in $X$, which is in bijection 
with the curve of B-lines in $F(Y)$ by \ref{-1}. 

Let $M_X(2;1,6)^1$ denote the locally closed subset 
of vector bundles that do not belong to 
$M_X(2;1,6)^0$. The complement of $M_X(2;1,6)^1$ in $M_X(2;1,6)$
is mapped bijectively by $c_2$ to  
the Fano surface $F(X)\subset J(X)$. 

Suppose that $M_X(2;1,6)^1$ is not empty. A fiber of $c_2$ passing 
through a point of $M_X(2;1,6)^1$ must contain another point, 
not in $M_X(2;1,6)^1$. Since the fibers are connected we deduce 
that $c_2$ contracts a curve in the moduli space to a same conic 
$q_{\ell}$. We thus have a five-dimensional family of elliptic
sextics $C$ mapping to $q_{\ell}$, that is, with a unisecant 
plane $\PP A$ such that $\sigma_{20}(\PP A)=C\cup D$, $D\subset 
G(2,H)$ for a hyperplane $H$,  and $G(2,H)\cap X=D\cup q_{\ell}$. 

But remember from the proof of Proposition \ref{-1} that the 
Palatini quartic $W$ meets the linear space $\PP^3_{\ell}$ along 
the union of two quadrics $Q_{\ell}$ and $Q^{\ell}$. The axis
of $D$ is contained in $W\cap\PP^3_{\ell}$, so we have a one
parameter family of such lines. The unisecant plane of $C$
contains the axis of  $D$, so we have three more parameters
for this plane 
and then the curve $C$ is determined. Hence a total of four 
parameters for the elliptic sextics mapping to a given conic
in $X$, and we deduce a contradiction. 

Since $M_X(2;1,6)^1=\emptyset$, 
 by \ref{vb-smooth} and \ref{m1-smooth} the moduli space
$M_X(2;1,6)$ is smooth.  Moreover $c_2$ defines
a bijection between the smooth surfaces $M_X(2;1,6)$ and $F(X)$, 
so it must be an isomorphism. \qed

\vspace{1cm}

{\small

\noindent 
{\bf Atanas Iliev}\\  
Institute of Mathematics, 
Bulgarian Academy of Sciences, 
Acad.G.Bonchev Str., bl.8\\  
1113 Sofia, Bulgaria\\
{\bf e-mail:} ailiev@math.bas.bg \\
Partially supported by grant MI-1503/2005
                  of the Bulgarian Foundation for Scientific Research

\bigskip

\noindent
{\bf Laurent Manivel}\\ 
Institut Fourier, 
Laboratoire de Math\'ematiques,
UMR 5582 (UJF-CNRS), BP 74\\ 
38402 St Martin d'H\`eres Cedex, France\\ 
{\bf e-mail:}  Laurent.Manivel@ujf-grenoble.fr
}

\end{document}